\documentclass[12pt]{article}

\usepackage[dvipsnames]{xcolor}
\usepackage{enumerate,enumitem}

\usepackage{amsmath,amssymb,amsfonts,mathrsfs, amsthm,mathtools, bm, bbm, dsfont, mathrsfs, amsthm}
\usepackage{graphicx, epstopdf}

\usepackage{fullpage}

\usepackage{ebgaramond}
\usepackage[OT1]{fontenc}
\usepackage[utf8]{inputenc}

\usepackage[colorlinks=true,breaklinks=true,bookmarks=true,urlcolor=MidnightBlue,citecolor=MidnightBlue,linkcolor=MidnightBlue,bookmarksopen=false,draft=false]{hyperref}
    
\newcommand{\LCON}{L_{{{\sf CO}}, N}}
\newcommand{\LEXN}{L_{{{\sf EX}}, N}}
\newcommand{\LCOSN}{L_{{{\sf CO, SYM}}, N}}
\newcommand{\LCO}{L_{{{\sf CO}}}}
\newcommand{\LEX}{L_{{{\sf EX}}}}
\newcommand{\LCOS}{L_{{{\sf CO, SYM}}}}
\newcommand{\PN}{$\mathcal{P}_{N}$}
\newcommand{\PIN}{$\mathcal{P}_{\infty}$}
\newcommand{\PNT}{$\mathcal{P}_{N}^{T}$}
\newcommand{\PIT}{$\mathcal{P}_{\infty}^{T}$}
\newcommand{\LPR}{L_{{{\sf PR}}}}
\newcommand{\LPRS}{L_{{{\sf PR, SYM}}}}

\newcommand{\R}{{\sf R}}

\newcommand{\N}{\mathcal{N}}
\newcommand{\M}{\mathcal{M}}

\allowdisplaybreaks

\usepackage{booktabs} 

\newtheorem{lemma}{Lemma}[]

\newtheorem{theorem}{Theorem}[]

\theoremstyle{definition}
\newtheorem{definition}{Definition}[]

\theoremstyle{assumption}
\newtheorem{assumption}{Assumption}[]

\theoremstyle{remark}

\begin{document}

\title{Nash Equilibria for Exchangeable Team against Team Games, their Mean Field Limit, and Role of Common Randomness}

\author{Sina Sanjari\thanks{Sina Sanjari is with the Department of Electrical and Computer Engineering at the University of Illinois Urbana-Champaign (UIUC). Email: sanjari@illinois.edu} \quad Naci Saldi\thanks{Naci Saldi is with the Department of Mathematics at Bilkent University, Turkey. Email:naci.saldi@bilkent.edu.tr} \quad  Serdar Y\"uksel\thanks{Serdar Y\"uksel is with the Department of Mathematics and Statistics at Queen's University, Canada. Email:yuksel@queensu.ca.}}

\date{}

\maketitle

\begin{abstract}
    We study stochastic mean-field games among finite number of teams with large finite as well as infinite number of decision makers. For this class of games within static and dynamic settings, we establish the existence of a Nash equilibrium, and show that a Nash equilibrium exhibits exchangeability in the finite decision maker regime and symmetry in the infinite one. To arrive at these existence and structural theorems, we endow the set of randomized policies with a suitable topology under various decentralized information structures, which leads to the desired convexity and compactness of the set of randomized policies. Then, we establish the existence of a randomized Nash equilibrium that is exchangeable (not necessarily symmetric) among decision makers within each team for a general class of exchangeable stochastic games. As the number of decision makers within each team goes to infinity (that is for the mean-field game among teams), using a de Finetti representation theorem, we show  existence of a randomized Nash equilibrium that is symmetric (i.e., identical) among decision makers within each team and also independently randomized. Finally, we establish that a Nash equilibrium for a class of mean-field games among teams (which is symmetric) constitutes an approximate Nash equilibrium for the corresponding pre-limit (exchangeable) game among teams with large but finite number of decision makers. We thus show that common randomness is not necessary for large team-against-team games, unlike the case with small sized teams. 
\end{abstract}

\section{Introduction.}

Stochastic teams entail a collection of decision makers (DMs) acting together to optimize a common cost function, but not necessarily sharing all the available information. At each time stage, each DM only has partial access to the global information which is defined as the \emph{information structure} (IS) of the team \cite{wit75}. If there is a pre-defined order according to which the DMs act, then the team is called a \emph{sequential team}. For sequential teams, if each DM's information depends only on primitive random variables, the team is \emph{static}. If at least one DM's information is affected by an action of another DM, the team is said to be \emph{dynamic}. If the cost or the probabilistic model are not equal or equivalent, such a setup is considered as a game. 
In this paper, we study a class of stochastic games where each player in the game entails a sequential team with  finite (but large) as well as countably infinite number of DMs. For this class of games among teams, we characterize existence and structural properties (i.e., symmetry) of a Nash equilibrium (NE). 


Mean-field games are limit models of symmetric non-zero-sum non-cooperative finite DM games with a mean-field interaction (see e.g., \cite{CainesMeanField2,CainesMeanField3,LyonsMeanField}). Mean-field games have many applications in financial engineering, economics, pricing in markets (see \cite{carmona2020applications} for a thorough review on some applications of mean-field games). The existence of a NE for mean-field games has been established in \cite{LyonsMeanField,bardi2019non,carmona2016mean,light2018mean,lacker2015mean,saldi2018markov}. Also, there have been several studies for mean-field games where the limits of sequences of Nash equilibria have been investigated as the number of DMs drives to infinity (see e.g., \cite{fischer2017connection,lacker2018convergence,bardi2014linear,LyonsMeanField,arapostathis2017solutions,lacker2016general}). In contrast to the setting of mean-field games has been studied in literature, in this paper, we study stochastic mean-field games among teams. Such problems exhibit additional nuances due to lack of convexity in team policies under decentralized ISs, see e.g., \cite{anantharam2007common,hogeboom2022zero,saldi2022geometry}. Such models have several applications as it combines mean-field game and team theory. 
A natural application area is in sensor network where a large collection of decentralized sensors (act as a large teams) shares their information (by their actions) to
a fusion center in the presence of jamming in the system entailing a large collection of decentralized jammers. As an another application, consider information sharing across a medium via interference channels involving a collection of multi-terminal encoder and decoder pairs where each encoder-decoder pair can be viewed as a team whose information serves as noise for the other teams. 

Unlike games where Nash equilibria is often of interest, the notion of optimality among each team is global optimality. In general, Nash equilibria (in the language of teams, person-by-person optimal solutions) only arrive at local optima and not global optima. This gap is specially significant for stochastic teams with countably infinite number of DMs (or mean-field teams). This is because, any perturbation of finitely many policies fails to significantly deviate the value of the expected cost. In view of this observation, the results on games may be inconclusive regarding global optimality for teams without uniqueness (e.g., see \cite{LyonsMeanField} for a strong sufficient condition of  monotonicity condition for establishing uniqueness). Non-uniqueness of equilibrium has been also studied; see e.g., \cite{bardi2019non,delarue2020selection,cardaliaguet2019example,hajek2019non,bayraktar2020non,cecchin2019convergence,lacker2018convergence}.


Mean-field teams under decentralized IS have been studied in \cite{sanjari2018optimal,sanjari2019optimal,SSYdefinetti2020}. For such models, existence and convergence of a globally optimal solution have been established in  \cite{sanjari2018optimal,sanjari2019optimal,SSYdefinetti2020}. In particular, it has been shown that under continuity of the cost function, mean-field team admit a symmetric (possibly randomized) globally optimal solution. Since we are studying games among teams, the NE concept needs to take into account global optimality notion among DMs within each team. In this paper, we study a class of mean-field games among finite number of teams with infinite DMs under the concept of NE which is globally optimal among DMs within each team but NE among teams (which are viewed as players). For this class of mean-field games, under sufficient continuity of the cost functions, we show that a NE exists that exhibits symmetry among DMs within teams. The existence of such a symmetric NE does not apriori follow directly from a fixed-point argument involving individual agents since restricting to symmetric policies in finding best response policies might be with a loss of global optimality for the DMs within each team. This is because, by fixing policies of other teams (players) to symmetric policies, each DM within a deviating team (player) faces a team problem that might not admit a symmetric optimal solution; e.g., see \cite[Example 1]{SSYdefinetti2020}. Related to our setting of games among large teams, in \cite{yu2021teamwise}, a mean-field type rank-based reward competition game has been studied among teams with continuum of players and centralized IS.  In addition, for zero-sum games involving finite teams against finite teams, the existence of a randomized saddle-point equilibrium with a common and private randomness has been established in \cite{hogeboom2022zero}. In this paper, our focus is on exchangeable games among teams with both large and infinite number of DMs to establish existence of a NE that is exchangeable for the finite case and symmetric in the mean-field limit. In \cite{hogeboom2022zero,anantharam2007common}, convexity of the set of randomized policies with a common and private randomness is required for the existence of an equilibrium, which can be realized by considering an arbitrary distribution over all possible common randomness variables. In practice, however, this may be too restrictive due to the arbitrary nature of distributions and since randomness must be externally provided.  

\subsection{Contributions.}

\begin{enumerate}
\item For general games among teams with finite number of DMs, in Theorem \ref{the:general}, we establish existence of a randomized NE with a common and independent randomness. For the special case of exchangeable games among teams with finite number of DMs, in Theorems \ref{the:1} and \ref{the:3} (for static and dynamic settings, respectively), we establish existence of a randomized NE, and show that it is exchangeable among DMs within each team (player). This NE might be asymmetric and also randomized, which can be correlated among DMs within each team (player).  Toward this goal, we show that the set of exchangeable randomized policies is convex and compact in Lemma  \ref{the:00}.

\item As the number of DMs within each team goes to infinity (that is for the mean-field game among teams), in Theorems \ref{the:2} and \ref{the:4} (for static and dynamic settings, respectively), we establish existence of a randomized NE, and show that it is symmetric (i.e., identical) and independently randomized among DMs within each team (player). We establish that the game problem can be relaxed to that where representative DMs' behavior can be regarded as the behavior of their corresponding symmetric population. This relaxes, in particular, the need for common randomness required for finite population regime.

\item In Theorems \ref{the:5} and \ref{the:6} (for static and dynamic settings, respectively), we show that a NE for a mean-field game among teams (which is symmetric and independent) constitutes an approximate NE for the corresponding game among teams with the mean-field interaction and large but finite number of DMs. 
\item As a final contribution, in view of our results noted above, our paper  presents a rigorous connection between Agent-Based-Modeling and Team-against-Team games, via the representative agents defining the game played in equilibrium. As noted in \cite[Section 4]{bertucci2019some}, Agent-based-Modeling is often used in economics, mathematical biology and other sciences to refer to a model in which macroscopic phenomena are captured by aggregating individual actions.
\end{enumerate}

\section{Problem formulation.}

\subsection{A generalized intrinsic model.}

Consider the class of stochastic games where DMs act in a pre-defined order (i.e., a sequential game \cite{wit75}). Under the intrinsic model for sequential games (described in discrete time), any action applied at any given time is regarded as applied by an individual DM, who acts only once. However, depending on the desired equilibrium concepts, IS, and cost functions, it is suitable to consider a collection of DMs as a single player acting as a team. In this setting, teams take part in a game. For this setting, Witsenhausen's intrinsic model is inadequate as it views each DM individually, which does not capture a joint deviation of a collection of DMs acting as a team, required for our equilibrium concept. To formalize this class of games, we introduce a class of $M$-player games using the generalized intrinsic model in \cite{SBSYgamesstaticreduction2021}, where each player consists of a collection of (one-shot) DMs either with finite or infinite members. A formal description has the following components:

\begin{itemize}[wide]
\item Let $\mathcal{M} :=\{1,2,\dots, M\}$ denote the set of players. For each $i\in \mathcal{M}$, let  $\N_{i}:=\{1,2,\dots, N_{i}\}$ denote the collection of DMs, DM$_{k}$ for $k\in \N_{i}$, acting as player $i$ (PL$^{i}$) (in other words, PL$^{i}$ encapsulates the collection of DMs indexed by $\mathcal{N}_{i}$). We denote  DM$_{k}$ of PL$^{i}$ by DM$^{i}_{k}$.

\item  There exists a collection of {\textit{measurable spaces}} $\{(\Omega, {\mathcal F}),(\mathbb{U}^i_{k},\mathcal{U}^i_{k}), (\mathbb{Y}^i_{k},\mathcal{Y}^i_{k})$, for $k\in \N_{i}$ and $i \in \mathcal{M}\}$, specifying the system's distinguishable events, DMs' control and observation spaces. The observation and action spaces are standard Borel spaces, and are described by the product spaces $\prod_{k \in \mathcal{N}^{i}}\mathbb{Y}^{i}_{k}$ and $\prod_{k \in \N_{i}}\mathbb{U}^{i}_{k}$, respectively for each PL$^{i}$.

\item  The $\mathbb{Y}^{i}_{k}$-valued observation variables are given by $y^{i}_{k} = h^{i}_{k} (\omega, \{{u}^{p}_{s}\}_{(s, p) \in L^{i}_{k}})$, where $L^{i}_{k}$ denotes the set of all DMs acting before DM$^{i}_{k}$ (i.e., $(s, p) \in  L^{i}_{k}$ if DM$^{p}_{s}$ acts before DM$^{i}_{k}$).

\item  An admissible policy for each PL$^{i}$ is denoted by $\underline{\gamma}^{i}:=\{\gamma^{i}_{k}\}_{k \in \mathcal{N}_{i}} \in {\Gamma}^{i}$ with $u_{k}^{i}=\gamma^{i}_{k}(y^{i}_{k})$ and the set of admissible policies for each PL$^{i}$ is described by the product space $\Gamma^{i}:=\prod_{k \in \mathcal{N}_{i}}\Gamma^{i}_{k}$, where $\Gamma^{i}_{k}$ is the set of all Borel measurable from $\mathbb{Y}^{i}_{k}$ to $\mathbb{U}^{i}_{k}$.

\item There is a {\textit{probability measure}} $\mathbb{P}$ on $(\Omega, {\cal F})$, making it a probability space. Let $\mathbb{E}$ denote the expectation with respect to $\mathbb{P}$. The prior probability measures in general can be subjective (for each DM or PL); however, in this paper, we will not discuss this, assuming that DMs have access to the common correct prior $\mathbb{P}$.
\end{itemize}

Under the above formulation, a sequential game is \emph{dynamic} if the
information available to at least one DM is affected by the action of at least one other DM. A game is \emph{static}, if for every DM the information is only affected by exogenous disturbances. ISs can also be categorized as \emph{classical}, \emph{quasi-classical} or \emph{non-classical}. An IS $\{y^i_{k}| i \in \M~~ \text{and}~~ k \in \N_{i}\}$ is classical if $y^i_{k}$ contains all of the information available to DM$^{p}_{s}$ who acts precedent to DM$^{i}_{k}$ (i.e., $(s, p) \in  L^{i}_{k}$). An IS is quasi-classical (or \emph{partially nested}), if whenever $u^{p}_{s}$  affects $y^i_{k}$ through the observation function $h^i_{k}$, $y^i_{k}$ contains $y^{p}_{s}$ (that is $\sigma(y^{p}_{s}) \subset \sigma(y^i_{k})$). An IS which is not partially nested is non-classical.

{In the paper, we will also allow for randomized policies, where in addition to $y^i_{k}$, each DM$^i_{k}$ has access to common and private randomization. This will be made precise in Section \ref{sec:strategic}. }

\subsection{Problems studied.}
We study stochastic games with each player involving finite (but large) number of DMs, or countably infinite number of DMs. For simplicity in our notations, we consider only two players in the game, denoted by PL$^{1}$ and PL$^{2}$ (i.e., $\M=\{1,2\}$). Our results remain valid for the the finite player setting. We address the following questions: 
\begin{enumerate}[label=\roman*),itemjoin={,\quad}]
   \item \emph{Does there exist a NE for games with finite number of DMs? Is this NE exchangeable?} We address the first question for general static games in Theorem \ref{the:general}, and the two questions in Theorem \ref{the:1} for exchangeable static games and in Theorem \ref{the:3} for exchangeable dynamic ones.  
    \item \emph{Does there exist a NE for games with countably infinite number of DMs? Does this NE admit symmetry properties?} We address these two questions in Theorem \ref{the:2} for mean-field static games and in Theorem \ref{the:4} for mean-field dynamic ones.
    
    \item \emph{Do Nash equilibria for games with countably infinite number of DMs constitute approximate Nash equilibria for the corresponding games with finite but large number of DMs?} We address this question in Theorem \ref{the:5} for static games and in Theorem \ref{the:6} for dynamic ones. 
\end{enumerate} 

The NE concept for such a class of games among teams should take into account the fact that DMs within players face a team problem with the desired global optimality notion for fixed policies of the other players. That is because, fixing policies of DMs within each player might lead to a local optimal and not global optimal. We additionally note that since the NE concept for such games among teams requires global optimality within teams, establishing the existence of a NE that exhibits exchangeability or symmetry requires further analysis that is not required for the classical stochastic games where each player is a singleton. This is because, in order to utilize a fixed point theorem, convexity and compactness of a set of policies (for each player, entailing all its DMs) under a topology are required. To address this difficulty, we endow a suitable topology on the set of randomized policies which leads to the convexity and compactness of a set of exchangeable randomized policies. As number of DMs goes to infinity, we use a de Finetti representation theorem and an argument used in \cite{SSYdefinetti2020} to establish symmetry, and to address the preceding questions.

\subsection{Static games.}
Our focus is on exchangeable stochastic games, and hence, we suppose that action and observation spaces are identical through DMs of each player and are subsets of appropriate dimensional Euclidean spaces, i.e., $\mathbb{U}^{i}_{k}=\mathbb{U}^{i}\subseteq\mathbb{R}^{n_{i}}$ and $\mathbb{Y}^{i}_{k}=\mathbb{Y}^{i}\subseteq\mathbb{R}^{m_{i}}$ for all $i \in \M$, where $n_{i}$ and $m_{i}$ are positive integers. In the following, we first introduce a class of exchangeable static games with finite number of DMs, and then, we introduce a class of static mean-field games.

\subsubsection{Finite DM static game \PN.}

Consider a stochastic game with finite number of DMs, i.e., for each PL$^{i}$,  ${\N_{i}}=\{1,\dots, N_{i}\}$,  $\underline{\gamma}_{N}^{i} := (\gamma^i_{1}, \cdots, \gamma_{N_{i}}^{i})$ and ${\Gamma}_{N}^{i} := \prod_{k=1}^{N_{i}} \Gamma^i_{k}$. Let the expected cost function for each PL$^{i}$ under a policy profile $\underline{\gamma}_{N}^{1:2}:=(\underline{\gamma}_{N}^{1}, \underline{\gamma}_{N}^{2})$ be given by
\begin{align}
J_{N}^{i}(\underline{\gamma}_{N}^{1:2}) &= \mathbb{E}^{\underline{\gamma}_{N}^{1:2}}\left[c^{i}(\omega_{0},\underline{u}_{N}^{1:2})\right]\label{eq:1.1}
\end{align}
for some Borel measurable cost function $c^{i}: \Omega_{0} \times \prod_{j=1}^{2} \prod_{k=1}^{N_{j}}\mathbb{U}^{j} \to \mathbb{R}_{+}$. We define $\omega_{0}$ as the $\Omega_{0}$-valued, cost function relevant, exogenous random variable, taking values from a Borel space $\Omega_{0}$ with its Borel $\sigma$-field $\mathcal{F}_{0}$. In the above, we used the notation $k:j$ to denote $\{k, \ldots, j\}$, and  $\underline {u}_{N}^{i}:={u}_{1:N_{i}}^{i}$, and $\mathbb{E}^{\underline{\gamma}_{N}^{1:2}}$ to denote the expectation with respect to $\mathbb{P}$ by replacing actions induced by policies $\underline{\gamma}_{N}^{1:2}$.

\begin{definition}[$\epsilon$-Nash equilibrium (NE) for \PN]
Let $\epsilon\geq 0$. A policy profile $\underline{\gamma}_{N}^{1*:2*}$ constitutes an $\epsilon$-NE if and only if the following inequalities hold for all $i\in \{1,2\}$  
\begin{flalign}\label{def:NE0}
J_{N}^{i}({\underline \gamma}^{1*:2*}_{N})\leq \inf_{\underline\gamma^{i}_{N}\in \Gamma^{i}_{N}}J_{N}^{i}({\underline \gamma}^{-i*}_{N}, \underline\gamma^{i}_{N}) +\epsilon,
\end{flalign} 
where $-i:=\{1,2\}\backslash \{i\}$. If $\epsilon=0$, the policy profile $\underline{\gamma}_{N}^{1*:2*}$ constitutes a NE. 
\end{definition}

We emphasize that in the above, DMs within a player optimize a common cost function but may have a decentralized IS. Hence, our notion of (player-wise) NE is a suitable equilibrium notion for such games since (global) optimality is desirable among DMs within players, which takes into account the fact that for each player, DMs face a team problem when policies of other players are fixed. In Section \ref{main:static}, we allow DMs to apply randomized policies, and hence, we re-write our game formulation and the equilibrium notion to incorporate randomized policies; see \eqref{eq:finitecost} (also \eqref{eq:pf}) and \eqref{eq:pinf}. 

For general exchangeable static games, the cost functions of the players satisfy the following exchangeability condition. 
\begin{assumption}\label{assump:exccost}
 $c^{i}$ is (separately) exchangeable with respect to actions of DMs for all $\omega_{0}$ and $i \in \M$, i.e., for any permutations $\sigma$ and $\tau$ of $\{1,\ldots,N_{1}\}$ and $\{1,\ldots,N_{2}\}$, respectively,  $$c^{i}\left(\omega_{0}, u^{1}_{1:N_{1}}, u^{2}_{1:N_{2}}\right)=c^{i}\left(\omega_{0}, u^{1}_{\sigma(1):\sigma(N_{1})}, u^{2}_{\tau(1):\tau(N_{2})}\right)\quad \forall \omega_{0},$$
 where $u^{1}_{\sigma(1):\sigma(N_{1})}=u^{1}_{\sigma(1)}, \ldots, u^{1}_{\sigma(N_{1})}$ and $u^{2}_{\tau(1):\tau(N_{2})}=u^{2}_{\tau(1)}, \ldots, u^{2}_{\tau(N_{2})}$.
\end{assumption}

 Let $\mathcal{P}(\cdot)$ denote the space of probability measures, and $\delta_{\{\cdot\}}$ denote the Dirac delta. An special case of the preceding game is when the costs satisfy Assumption \ref{assump:exccost} and are given by
\begin{equation}\label{eq:exceee}
c^{i}\left(\omega_{0}, u^{1}_{1:N_{1}}, u^{2}_{1:N_{2}}\right)=\frac{1}{N_{i}}\sum_{k=1}^{N_{i}}c^{i}_{k}\left(\omega_{0},u^{i}_{k},\Xi^1(\frac{1}{N_{1}}\sum_{p=1}^{{N}_{1}}\delta_{u^{1}_{p}}), \Xi^2(\frac{1}{N_{2}}\sum_{p=1}^{{N}_{2}}\delta_{u^{2}_{p}})\right),
\end{equation}
for some weakly continuous functions $\Xi^{i}:\mathcal{P}(\mathbb{U}^{i}) \to \tilde{\mathbb{U}}^{i}$ with $\tilde{\mathbb{U}}^{i}$, which is a subset of appropriate dimensional Euclidean space, and for some Borel measurable cost functions $c^{i}_{k}:\Omega_{0} \times \mathbb{U}^{i} \times \prod_{j=1}^{2} \tilde{\mathbb{U}}^{j}  \rightarrow \mathbb{R}_{+}$ that satisfies
\begin{flalign}
&c^{i}_{k}\left(\omega_{0},u^{i}_{k},\Xi^1(\frac{1}{N_{1}}\sum_{p=1}^{{N}_{1}}\delta_{u^{1}_{p}}), \Xi^2(\frac{1}{N_{2}}\sum_{p=1}^{{N}_{2}}\delta_{u^{2}_{p}})\right)=c^{i}_{\beta(k)}\left(\omega_{0},u^{i}_{\beta(k)},\Xi^1(\frac{1}{N_{1}}\sum_{p=1}^{{N}_{1}}\delta_{u^{1}_{p}}), \Xi^2(\frac{1}{N_{2}}\sum_{p=1}^{{N}_{2}}\delta_{u^{2}_{p}})\right),\label{costexmf}
\end{flalign}
for every permutation $\beta=\sigma$ or $\beta=\tau$.
Clearly, if $c^{i}_{k}=\tilde{c}^{i}$, then the above exchangeability condition holds. When the number of DMs drives to infinity, we consider the special case of $c^{i}_{k}=\tilde{c}^{i}$. We note that weak continuity of $\Xi^{i}$ is not required for exchangeability condition, but it is required for our main theorems.

\subsubsection{Mean-field static game \PIN.}

Consider a stochastic game with countably infinite number of DMs, i.e., let $\mathcal{N}_{i}=\mathbb{N}$, ${\Gamma}^{i}:=\prod_{k \in \mathbb{N}} \Gamma^{i}_{k}$ and $\underline{\gamma}^{i}:=(\gamma^{i}_{1},\gamma^{i}_{2},\ldots)$. Let the expected cost of PL$^{i}$ under a policy profile $\underline{\gamma}^{1:2}$ be given by
\begin{equation}\label{eq:2.5.5}
J^{i}_{\infty}(\underline{\gamma}^{1:2})=\limsup\limits_{N_{1}, N_{2}\rightarrow \infty}  \mathbb{E}^{\underline{\gamma}^{1:2}}\left[\frac{1}{N_{i}}\sum_{k=1}^{N_{i}}c^{i}\left(\omega_{0},u^{i}_{k},\Xi^{1}(\frac{1}{N^{1}}\sum_{p=1}^{{N}_{1}}\delta_{u^{1}_{p}}),\Xi^{2}( \frac{1}{N_{2}}\sum_{p=1}^{{N}_{2}}\delta_{u^{2}_{p}})\right)\right],
\end{equation}
for some weakly continuous functions $\Xi^{i}:\mathcal{P}(\mathbb{U}^{i}) \to \tilde{\mathbb{U}}^{i}$ and Borel measurable cost functions $c^{i}_k:\Omega_{0} \times \mathbb{U}^{i} \times \prod_{j=1}^{2} \tilde{\mathbb{U}}^{j}  \rightarrow \mathbb{R}_{+}$. We assume that $N_{1}$ and $N_{2}$ converge with at same rate to infinity, i.e., there exist constants $a,b\in \mathbb{R}$ such that $a\leq \liminf\limits_{N_{1}, N_{2}\rightarrow \infty} \frac{N_{1}}{N_{2}}\leq \limsup\limits_{N_{1}, N_{2}\rightarrow \infty} \frac{N_{1}}{N_{2}}\leq b.$

\begin{definition}[Nash equilibrium (NE) for \PIN]
A policy profile $\underline{\gamma}^{1*:2*}$ is NE (or mean-field equilibrium) for \PIN\ if and only if the following inequalities hold for all $i\in \{1,2\}$  
\begin{flalign}\label{def:NE}
J^{i}_{\infty}({\underline \gamma}^{1*:2*})\leq J^{i}_{\infty}({\underline \gamma}^{-i, *}, \underline\gamma^{i})\qquad \forall \underline\gamma^{i}\in \Gamma^{i}.
\end{flalign}
\end{definition}

Our first main goal is to establish the existence of a symmetric (identical) randomized NE for static mean-field games ${\mathcal{P}_{\infty}}$.  
{Our results require the following absolute continuity condition. 

\begin{assumption}\label{assump:ind}
For every $N_{1}, N_{2}\in \mathbb{N}\cup \{\infty\}$, let $\tilde{\mu}^{N}$ be the conditional distribution of observations $\underline{y}^{1:2}$ given $\omega_{0}$. There exists a probability measure $Q^{i}_{k}$ on $\mathbb{Y}^{i}$ and a bounded function $f^{i}_{k}:\mathbb{Y}^{i}\times \Omega_{0}\times \prod_{j=1,\not=k}^{N_{i}}\mathbb{Y}^{i} \times \prod_{j=1}^{N_{-i}}\mathbb{Y}^{-i}\to \mathbb{R}_{+}$ for all $i\in \M$ and $k\in \N_{i}$ such that for all Borel set $B^{i}_{k}$ in $\mathbb{Y}^{i}$ (with $B :=\prod_{i=1}^{2} B^{i}_{1}\times \cdots \times B^{i}_{N_{i}}$)
\begin{flalign}
&\tilde{\mu}^{N}(B \big|\omega_{0})=\prod_{i=1}^{2}\prod_{k=1}^{N_{i}}\int_{B^i_{k}} f^{i}_{k}\left(y^{i}_{k}, \omega_{0},y^{i}_{-k}, \underline{y}^{-i}\right)Q^{i}_{k}(dy^{i}_{k})\label{eq:abscon},
\end{flalign}
where $y^{i}_{-k}=\{y^{i}_{1:N_{i}}\}\backslash \{y^{i}_{k}\}$.
 \end{assumption}

Under this assumption, via change of measure argument (see e.g., \cite{wit88,SBSYgamesstaticreduction2021}), we can equivalently view the observations of each DM as independent and also independent of $\omega_{0}$. The above allows us to introduce a suitable topology under which the space of randomized policies is Borel (see Section \ref{sec:strategic}). In addition, our main results (i.e., Theorems \ref{the:1} and \ref{the:2}) for static games impose the following assumptions on the observations, action spaces, and cost functions. 

\begin{assumption}\label{assump:oc}
\hfill 
\begin{itemize}
\item [(i)] $(y^{1}_{k})_{k\in \N_{1}}$ and $(y^{2}_{k})_{k\in \N_{2}}$ are independent, conditioned on $\omega_{0}$;
\item [(ii)] For all $i \in \M$, $(y^{i}_{k})_{k\in \N_{i}}$ have an identical distribution, conditioned on $\omega_{0}$.
\end{itemize}
\end{assumption}

\begin{assumption}\label{assump:cc}
For $i\in \M$,
\begin{itemize}
\item [(i)] $\mathbb{U}^{i}$ is compact.
\item [(ii)] $c^{i}(\omega_{0}, \cdot, \cdot, \cdot)$ in \eqref{eq:exceee}, is continuous for all $\omega_{0}$.
\end{itemize}
\end{assumption}

We note that under Assumption \ref{assump:oc}, we can rewrite \eqref{eq:abscon} as follows:
\begin{flalign}
&\tilde{\mu}^{N}(B \big|\omega_{0})=\prod_{i=1}^{2}\prod_{k=1}^{N_{i}}\int_{B^i_{k}} \hat{f}^{i}\left(y^{i}_{k}, \omega_{0}\right)Q^{i}(dy^{i}_{k})\label{eq:abscon-iden},
\end{flalign}
where $\hat{f}^{i}:\mathbb{Y}^{i}\times \Omega_{0} \to \mathbb{R}_{+}$ and $Q^{i}$ are identical for all DMs within player $i$ for $i=1,2$.

In Section \ref{main:static}, we establish the existence of a randomized NE for \PIN, and we show that this NE is symmetric (identical among DMs within players). To this end, we first establish existence of a randomized NE for \PN\ that is exchangeable, and then, we use this result to arrive at existence result for \PIN. For our results in Section \ref{main:static}, we let Assumptions \ref{assump:exccost} and  \ref{assump:ind} to hold.} For exchangeable games with finite number of DMs, we consider Assumption \ref{assump:cc}, but we relax Assumption \ref{assump:oc}, where we allow DMs to have   (finite)-exchangeable  that can be correlated (see Assumption \ref{assump:obsx}). When number of DMs is driven to infinity, we consider mean-field interaction among them, and we let Assumptions \ref{assump:oc} and  \ref{assump:cc} to hold.

\subsection{Dynamic games.}\label{sec:dn}

Again, our focus is on exchangeable games but dynamic ones. We let action, observation, and state spaces, respectively, be identical through DMs $k\in \mathcal{N}_{i}$, and for simplicity, also identical over time $t=0,\dots,T-1$, and subset of appropriate dimensional Euclidean spaces $\mathbb{U}^{i}_{k, t}=\mathbb{U}^{i}\subseteq\mathbb{R}^{n}$, $\mathbb{Y}^{i}_{k,t}=\mathbb{Y}^{i}\subseteq\mathbb{R}^{n^{\prime}}$, $\mathbb{X}^{i}_{k, t}=\mathbb{X}^{i}\subseteq\mathbb{R}^{n^{\prime\prime}}$, where $n$, $n^{\prime}$ and $n^{\prime\prime}$  are positive integers. In the following, we first introduce a class of games with finite number of DMs, and then, we introduce a class of dynamic mean-field games. State dynamics and observations of DMs are given by
\begin{flalign}
x_{k, t+1}^{i}&=f_{t}^{i}\bigg(x_{k, t}^{i},u_{k, t}^{i},\Xi^{1}_{x}(\frac{1}{N_{1}}\sum_{p=1}^{N_{1}}\delta_{x_{p, t}^{1}}), \Xi^{2}_{x}(\frac{1}{N_{2}}\sum_{p=1}^{N_{2}}\delta_{x_{p, t}^{2}}),\Xi^{1}_{u}(\frac{1}{N_{1}}\sum_{p=1}^{N_{1}}\delta_{u_{p, t}^{1}}), \Xi^{2}_{u}(\frac{1}{N_{2}}\sum_{p=1}^{N_{2}}\delta_{u_{p, t}^{2}}), w_{k, t}^{i}\bigg),\label{eq:mfdynamics2}\\
y_{k, t}^{i}&=h_{t}^{i}\left(x_{k, 0:t}^{i},u_{k, 0:t-1}^{i},v_{k, 0:t}^{i}\right)\label{eq:mfobs2},
\end{flalign}
for some weakly continuous functions $\Xi^{i}_{x}:\mathcal{P}(\mathbb{X}^{i})\to \tilde{\mathbb{X}}^{i}$ and $\Xi^{i}_{u}:\mathcal{P}(\mathbb{U}^{i})\to \tilde{\mathbb{U}}^{i}$, where $\tilde{\mathbb{X}}^{i}$ and $\tilde{\mathbb{U}}^{i}$ are subsets of appropriate dimensional Euclidean spaces. In the above, functions $f_{t}^{i}$ and $h_{t}^{i}$ are measurable functions and $v_{k,t}^{i}$ and $w_{k,t}^{i}$ correspond to uncertainties in state dynamics and observations. 

\subsubsection{Finite DM dynamic game \PNT.}

Consider a stochastic dynamic game with finite number of DMs. Each PL$^{i}$'s admissible policies $\pmb{\underline\gamma^{i}_{N}}:=\underline\gamma^{i}_{N, 0:T-1}$ with $\underline\gamma^{i}_{N, t}:=\gamma^{i}_{1:N_{i},t}$ are measurable functions so that $u^i_{k, t} = \gamma^i_{k, t}(y^i_{k, t})$ for all $k\in \mathcal{N}_{i}$ and $t=0,\dots, T-1$. We use bold letters for variables over time $t=0,\dots, T-1$, and underline variable to denote DMs within a player $k=1, \ldots, N_{i}$. Denote the space of policies for each player with $\pmb{\Gamma^{i}_{N}}=\prod_{t=0}^{T-1}\prod_{k=1}^{N_{i}}\Gamma^{i}_{k,t}$. Let the expected cost function of $\pmb{\underline \gamma^{1:2}_{N}}$ for \PNT\ be given by
  \begin{flalign}
J_{N}^{i, T}(\pmb{\underline \gamma^{1:2}_{N}})&=\mathbb{E}^{\pmb{\underline \gamma^{1:2}_{N}}}\bigg[\frac{1}{N_{i}}\sum_{k=1}^{N_{i}}\sum_{t=0}^{T-1}c^{i}_{k}\bigg(\omega_{0},x_{k, t}^{i},u_{k, t}^{i}, \Xi^{1}_{x}(\frac{1}{N_{1}}\sum_{p=1}^{N_{1}}\delta_{x_{p, t}^{1}}),\label{eq:mfcost}\\
&\qquad \qquad\Xi^{2}_{x}(\frac{1}{N_{2}}\sum_{p=1}^{N_{2}}\delta_{x_{p, t}^{2}}), \Xi^{1}_{u}(\frac{1}{N_{1}}\sum_{p=1}^{N_{1}}\delta_{u_{p, t}^{1}}), \Xi^{2}_{u}(\frac{1}{N_{2}}\sum_{p=1}^{N_{2}}\delta_{u_{p, t}^{2}})\bigg)\bigg]\nonumber,
\end{flalign}
for some weakly continuous functions $\Xi^{i}_{x}:\mathcal{P}(\mathbb{X}^{i})\to \tilde{\mathbb{X}}^{i}$ and $\Xi^{i}_{u}:\mathcal{P}(\mathbb{U}^{i})\to \tilde{\mathbb{U}}^{i}$, where $\tilde{\mathbb{X}}^{i}$ and $\tilde{\mathbb{U}}^{i}$, and for some Borel measurable cost functions $c_{k}^{i}:\Omega_{0}\times \mathbb{X}^{i}\times \mathbb{U}^{i}\times \prod_{j=1}^{2}\tilde{\mathbb{X}}^{j}\times \tilde{\mathbb{U}}^{j}$ that satisfies
\begin{flalign}
c^{i}_{k}(\cdot,x_{k, t}^{i},u_{k, t}^{i}, \cdot, \ldots, \cdot) = c^{i}_{\beta(k)}(\cdot, x_{\beta(k),t}^{i}, u_{\beta(k),t}^{i}, \cdot, \ldots, \cdot) \label{eq:costexdymf}
\end{flalign}
for permutations $\beta=\sigma$ or $\beta=\tau$. Random variable $\omega_{0}$ is again a cost-related random variable. If $c^{i}_{k}=\tilde{c}^{i}$, then \eqref{eq:costexdymf} holds. When the number of DMs drives to infinity, we consider the special case $c^{i}_{k}=\tilde{c}^{i}$.

\begin{definition}[$\epsilon$-NE for \PNT]
Let $\epsilon\geq 0$. A policy profile $\pmb{\underline \gamma^{1*:2*}_{N}}$ constitutes an $\epsilon$-NE if and only if the following inequalities hold for all $i\in \{1,2\}$  
\begin{flalign}\label{def:NE1}
J_{N}^{i, T}(\pmb{\underline \gamma^{1*:2*}_{N}})\leq \inf_{\pmb{\underline\gamma^{i}_{N}}\in \pmb{\Gamma^{i}_{N}}} J_{N}^{i,T}(\pmb{\underline\gamma^{-i*}_{N}}, \pmb{\underline\gamma^{i}_{N}})+\epsilon.
\end{flalign}  
If $\epsilon=0$, then the policy profile $\pmb{\underline \gamma^{1*:2*}_N}$ constitutes a NE.
\end{definition}

 In Section \ref{sec:4}, we allow DMs to apply randomized policies, and re-write our game formulation and the equilibrium notion to incorporate randomized policies; see \eqref{eq:finiterandom} and \eqref{eq:infiniterandom}. 

\subsubsection{Mean-field dynamic game \PIT.}

Consider a stochastic dynamic game with countably infinite number of DMs. Let $\mathcal{N}_{i}=\mathbb{N}$, $\pmb{\Gamma^{i}}=\prod_{t=0}^{T-1}\prod_{k\in \mathbb{N}}\Gamma^{i}_{k,t}$ and $\pmb{\underline\gamma^{i}}:=\underline\gamma^{i}_{0:T-1}$ with $\underline\gamma^{i}_{t}:=(\gamma^{i}_{1,t}, \gamma^{i}_{2,t}, \ldots)$. Let the expected costs of a policy profile $\pmb{\underline \gamma^{1:2}}$ for \PIT\ be given by
 \begin{flalign}
&J_{\infty}^{i,T}(\pmb{\underline \gamma^{1:2}})=\limsup\limits_{N_{1},N_{2}\rightarrow \infty}J_{N}^{i,T}(\pmb{\underline \gamma^{1:2}_{N}})\label{eq:dmfcost}.
\end{flalign}

 Similar to \eqref{def:NE1}, we can define NE (or mean-field equilibrium) for \PIT. Our solution approach for the dynamic setting is similar to that of static one with additional technical arguments and assumptions. Our theorems for the dynamic games (Theorems \ref{the:3} and \ref{the:4}) require that absolute continuity conditions in Assumption \ref{assump:ind1} hold. Assumption \ref{assump:ind1} allows us to endow a suitable topology for the space of randomized policies. Furthermore, Theorem \ref{the:4} requires the following assumption.

\begin{assumption}\label{assump:c}
For $i=1,2$,
\begin{itemize}
\item[(i)] $f_{t}^{i}(\cdot, \ldots, \cdot, w_{k,t}^{i})$ and $h_{t}^{i}(\cdot,\cdot, v_{k,0:t}^{i})$ are continuous and uniformly bounded for all $(w_{k,t}^{i}, v_{k,0:t}^{i})$ and for $t=0,\dots,T-1$.
\item [(ii)] $c^{i}(\omega_{0}, \cdot, \ldots, \cdot)$ is continuous and uniformly bounded for all $\omega_{0}$.
\item [(iii)] $\mathbb{U}^{i}$ is compact.
\end{itemize}
 \end{assumption}
\begin{assumption}\label{assump:2}
For $i=1,2$,
\begin{itemize}
\item[(i)] $(x_{k,0}^{i})_{k\in \mathcal{N}_{i}}$ are i.i.d., conditioned on $\omega_{0}$;  
\item[(ii)] $(w^{i}_{k, t})_{k\in \mathcal{N}_{i}}$ are i.i.d. for  $t=0,\dots,T-1$, and $(w^{i}_{k, 0:T-1})$ are mutually independent, and independent of $\omega_{0}$ and $(x_{k,0}^{i})_{k\in \mathcal{N}_{i}}$. $(v^{i}_{k, t})_{k\in \mathcal{N}_{i}}$ are i.i.d. for  $t=0,\dots,T-1$, and $(v^{i}_{k, 0:T-1})$ are mutually independent, and independent of $\omega_{0}$, $(x_{k,0}^{i})_{k\in \mathcal{N}_{i}}$. and $(w^{i}_{k, t})_{k\in \mathcal{N}_{i}}$ for  $t=0,\dots,T-1$.
\end{itemize}
\end{assumption}

 In Section \ref{sec:4}, we establish the existence of a  randomized NE for \PIT, and we show that it is symmetric (identical among DMs within players). Similar to the static case, for a more general formulation of \PNT, we first establish existence of a randomized NE that is exchangeable, and then, we use the result to arrive at existence result for \PIT. For our results in Section \ref{sec:4}, we let Assumption \ref{assump:ind1} to hold, but we let Assumptions \ref{assump:ind2}, \ref{assump:c}, and \ref{assump:2} hold only when they are needed.

\section{Exchangeable Static Games.}\label{main:static}
In this section, we study exchangeable static stochastic games.  We first introduce randomized policies with their suitable topology, and then, we establish existence of a NE for \PN\ that is exchangeable, and a NE for \PIN\ that is symmetric.

\subsection{Topology on control policies for static games.}\label{sec:strategic}
In the following, we introduce randomized policies as Borel probability measures, equipped with a suitable topology. Following \cite{yuksel2018general,wit88}, via a change of measure argument, Assumption \ref{assump:ind} enables us to equivalently view observations of DMs of players independent and also independent of $\omega_{0}$. Under Assumption \ref{assump:ind}, we  separate DMs' policy spaces (both across the player and DMs), equip them with a suitable topology. Let  
 \begin{flalign}
 \Theta^{i}_{k}:=&\bigg\{P \in \mathcal{P}(\mathbb{U}^{i}\times \mathbb{Y}^{i})\bigg|\label{eq:topogamd}\\
 &P(B)=\int_{B}\delta_{\{g^{i}_{k}(y^{i}_{k})\}}(du^{i}_{k})Q^{i}_{k}(dy^{i}_{k}),~g^{i}_{k}:\mathbb{Y}^{i}\to \mathbb{U}^{i}, ~ B\in \mathcal{B}(\mathbb{U}^{i}\times \mathbb{Y}^{i})\bigg\}\nonumber.
 \end{flalign}
Also, let
\begin{flalign}
\mathcal{R}^{i}_{k}:=\left\{P \in \mathcal{P}(\mathbb{U}^{i}\times \mathbb{Y}^{i}) \bigg| P(B)=\int_{B}\Pi^{i}_{k}(du^{i}_{k}|y^{i})Q^{i}_{k}(dy^{i}_{k}), ~ B\in \mathcal{B}(\mathbb{U}^{i}\times \mathbb{Y}^{i})\right\}\label{eq:topogam},
\end{flalign}
where $\Pi^{i}_{k}$ is a stochastic kernel from $\mathbb{Y}^{i}$ to $\mathbb{U}^{i}$. The set $\Theta^{i}_{k}$ is the set of extreme points of $\mathcal{R}^{i}_{k}$. Hence, $\Theta^{i}_{k}$ inherits Borel measurability and topological properties of the Borel measurable set $\mathcal{R}^{i}_{k}$ under the weak convergence topology \cite{BorkarRealization}. We identify the set of relaxed policies\footnote{This set corresponds to the set of Young measures in \cite{young1937generalized}.} for each DM$^{i}_{k}$, $\Gamma^{i}_{k}$, by $\mathcal{R}^{i}_{k}$ with the convergence in policies is defined by
\begin{flalign}\label{eq:convergence}
\gamma^{i}_{n, k} \xrightarrow[n\to \infty]{}  \gamma^{i}_{k} \iff \gamma^{i}_{k, n}(du^{i}_{k}|y^{i}_{k})Q^{i}_{k}(dy^{i}_{k}) \xrightarrow[n\to \infty]{\text{weakly}} \gamma^{i}_{k}(du^{i}_{k}|y^{i}_{k})Q^{i}_{k}(dy^{i}_{k}).
\end{flalign}
We equip $\mathcal{P}(\mathbb{Y}^{i}\times \mathbb{U}^{i})$ with the $w$-$s$ topology which is the coarsest topology under which the map $\int \kappa(y^i_{k}, u^i_{k})P(dy^i_{k}, du^i_{k}): \mathcal{P}(\mathbb{Y}^{i}\times \mathbb{U}^{i}) \to \mathbb{R}$ is continuous for every measurable and bounded function $\kappa$, which is continuous in $u^i_{k}$ for every $y^i_{k}$. Unlike the weak convergence topology, $\kappa$ need not to be continuous in $y^i_{k}$ (see \cite{Schal}). Since the marginals on observations are fixed, under the $w$-$s$ topology, the convergence coincides with that in the weak convergence topology, i.e., the convergence of probability measures is the weak convergence (see \cite[Theorem 3.10]{Schal}). Hence. we can view convergence in \eqref{eq:convergence} in terms of $w$-$s$ topology without any loss of generality.

\subsubsection{Randomized policies for \PN.}    

 The above formulation of relaxed policies for $\Gamma^{i}_{N}:=\prod_{k=1}^{N_{i}}\Gamma^{i}_{k}$ allows us to introduce the set randomized policies $L_{N}^{i}:=\mathcal{P}(\Gamma_{N}^{i})$ for each PL$^i$  as a collection of Borel probability measures on $\Gamma_{N}^{i}$, where Borel $\sigma$-field $\mathcal{B}(\Gamma^{i}_{k})$ is induced by the introduced topology in \eqref{eq:convergence}. In the following, we introduce different sets of randomized policies by allowing common and individual randomness among DMs. Let 
\begin{flalign*}
\LCON^{i}:=\bigg\{P_{\pi} \in L_{N}^{i}&\bigg{|}\text{for all}~A_{k} \in \mathcal{B}(\Gamma^{i}_{k}):P_{\pi}(\gamma^{i}_{1} \in A_{1},\dots,\gamma^{i}_{N_{i}}\in A_{N_{i}}) \nonumber\\
&=\int_{z^{i}\in [0,1]}\prod_{k=1}^{N_{i}}P_{\pi,  k}(\gamma^{i}_{k}\in A_{k}|z^{i})\eta^{i}(dz^{i}), ~~~~\eta^{i} \in \mathcal{P}([0, 1])\bigg\},
\end{flalign*} 
where $\eta^{i}$ is the distribution of common randomness (independent from the intrinsic exogenous system random variables). In the above, for every fixed $z^{i}$, $P_{\pi, k}\in \mathcal{P}(\Gamma^{i}_{s})$ corresponds to an independent randomized policy for each DM$^{i}_{k}$ ($k\in \N_{i}$ and $i\in \M$). Since $\LCON^{i}$ and $L_{N}^{i}$ are equal (see \cite[Theorem A.1]{SSYdefinetti2020}), $L_{N}^{i}$ corresponds to randomized policies induced by individual and common randomness. 
{We next recall the definition of {\it exchangeability} for random variables.
\begin{definition}
Random variables $x^{1:m}$ defined on a common probability space are $m$-\it{exchangeable} if for any permutation $\sigma$ of the set $\{1,\dots,m\}$, $\mathcal{L}\left(x^{\sigma(1):\sigma(m)}\right)=\mathcal{L}\left(x^{1:m}\right)$, where $\mathcal{L}(\cdot)$ denotes the (joint) law of random variables. Random variables $(x^{1},x^{2},\dots)$ are {\it infinitely-exchangeable} if finite distributions of $(x^{1},x^{2},\dots)$ and $(x^{\sigma(1)},x^{\sigma(2)},\dots)$ are identical for any finite permutation (affecting only finitely many elements) of $\mathbb{N}$. 
\end{definition}}

The set of exchangeable randomized policies $\LEXN^{i}\subseteq L_{N}^{i}$ for each PL$^i$ is given by
\begin{flalign}
\LEXN^{i}:=\bigg\{&P_{\pi} \in L_{N}^{i}\bigg{|} \forall A_{k} \in \mathcal{B}(\Gamma^{i}_{k})~\text{and}~\forall\sigma \in S_{N_{i}}:\nonumber\\
&P_{\pi}(\gamma^{i}_{1} \in A_{1},\dots,\gamma^{i}_{N_{i}}\in A_{N_{i}})=P_{\pi}(\gamma_{\sigma(1)}^{i} \in A_{1},\ldots,\gamma_{\sigma(N_{i})}^{i}\in A_{N_{i}})\bigg\}\label{eq:LEXN},
\end{flalign}
where $S_{N_{i}}$ is the set of all permutations of $\{1,\ldots,N_{i}\}$. We note that $\LEXN^{i}$ is a convex subset of $L^{i}_{N}$.  {Let  $\LCOSN^{i}$ be the set of all symmetric (identical) randomized policies induced by a common and individual randomness:
\begin{flalign*}
\LCOSN^{i}:=\bigg\{P_{\pi} \in L^{i}_{N}&\bigg{|}\forall A_{k} \in \mathcal{B}(\Gamma^{i}_{k}):P_{\pi}(\gamma^{i}_{1} \in A_{1},\dots,\gamma^{i}_{N_{i}}\in A_{N_{i}}) \nonumber\\
&=\int_{z^{i}\in [0,1]}\prod_{k=1}^{N_{i}}\tilde{P}_{\pi}(\gamma^{i}_{k} \in A_{k}|z)\eta^{i}(dz^{i}), ~~~~\eta^{i} \in \mathcal{P}([0, 1])\bigg\},
\end{flalign*}
where for all $k \in \N_{i}$, conditioned on $z^{i}$, $\tilde{P}_{\pi} \in \mathcal{P}(\Gamma^{i}_{k})$ corresponds to an independent randomized policy for each DM$^{i}_{k}$ ($k \in \N_{i}$), which is symmetric among PL$^{i}$'s DMs. Note that $\LCOSN^i \subset \LEXN^i$. 

\subsubsection{Randomized policies for \PIN.}    

When the number of DMs for each PL is countably infinite, we define the corresponding sets of randomized policies $L^{i}, \LCO^{i}, \LEX^{i}$, and $\LCOS^{i}$ using Ionescu Tulcea
extension theorem and by iteratively adding new coordinates for our probability measures (see e.g., \cite{InfiniteDimensionalAnalysis,HernandezLermaMCP}). The set of randomized policies $L$ on the infinite product Borel spaces $\Gamma^{i}:=\prod_{k\in \mathbb{N}}\Gamma^{i}_{k}$ is defined as
$L^{i}:=\mathcal{P}(\Gamma^{i})$.
Let the set of all randomized policies with a common and independent randomness be given by
\begin{flalign*}
\LCO^{i}:=\bigg\{P_{\pi} \in L^{i}&\bigg{|}\text{for all}~A_{k} \in \mathcal{B}(\Gamma^{i}_{k}): P_{\pi}(\gamma^{i}_{1} \in A_{1},\gamma^{i}_{2} \in A_{2},\dots)\nonumber\\
&=\int_{z^{i}\in [0,1]}\prod_{k\in \mathbb{N}}P_{\pi, k}(\gamma^{i}_{k}\in A_{k}|z^{i})\eta^{i}(dz^{i}), ~~~~\eta^{i} \in \mathcal{P}([0, 1])\bigg\}.
\end{flalign*}
Let the set of all infinitely exchangeable randomized policies $\LEX^{i}$ be given by
\begin{flalign*}
\LEX^{i}:=\bigg\{&P_{\pi} \in L^{i}\bigg{|}\text{for all}~A_{k} \in \mathcal{B}(\Gamma^{i}_{k})~\text{and for all $N_{i}\in\mathbb{N}$, and for all $\sigma \in S_{N_{i}}$:} \nonumber\\
&P_{\pi}(\gamma^{i}_{1} \in A_{1},\dots,\gamma^{i}_{N_{i}}\in A_{N_{i}})=P_{\pi}(\gamma^{i}_{\sigma(1)} \in A_{1},\dots,\gamma^{i}_{\sigma(N_{i})}\in A_{N_{i}})\bigg\}.
\end{flalign*}
Let the symmetric policies with a common and independent randomness be given by
\begin{flalign*}
\LCOS^{i}:=\bigg\{P_{\pi} \in L^{i}&\bigg{|}\text{for all}~A_{k} \in \mathcal{B}(\Gamma^{i}_{k}): P_{\pi}(\gamma^{i}_{1} \in A_{1},\gamma^{i}_{2}\in A_{2}, \dots)\nonumber\\
&=\int_{z^{i}\in [0,1]}\prod_{k\in \mathbb{N}}\tilde{P}_{\pi}(\gamma^{i}_{k} \in A_{k}|z^{i})\eta^{i}(dz^{i}), ~~~~\eta^{i} \in \mathcal{P}([0, 1])\bigg\}.
\end{flalign*}

For the infinite population setting, we also define the following set of randomized policies with only private independent randomness:
\begin{flalign}
\LPR^{i}:=\bigg\{&P_{\pi} \in L^{i}\bigg{|}\text{for all}~A_{k} \in \mathcal{B}(\Gamma^{i}_{k}):\\
&P_{\pi}(\gamma^{i}_{1} \in A_{1},\gamma^{i}_{2} \in A_{2},\dots)=\prod_{k\in \mathbb{N}}P_{\pi, k}(\gamma^{i}_{k}\in A_{k}),~\text{for}~P_{\pi, k}\in \mathcal{P}(\Gamma^{i}_{k})\bigg\}\nonumber.
\end{flalign} 
Define the set of symmetric randomized policies with private independent randomness as
\begin{flalign*}
\LPRS^{i}:=\bigg\{&P_{\pi} \in L^{i}\bigg{|}\text{for all}~A_{k} \in \mathcal{B}(\Gamma^{i}_{k}):\\
&P_{\pi}(\gamma^{i}_{1} \in A_{1},\gamma^{i}_{2}\in A_{2}, \dots)=\prod_{k\in \mathbb{N}}\tilde{P}_{\pi}(\gamma^{i}_{k} \in A_{k}),~\text{for}~\tilde{P}_{\pi} \in \mathcal{P}(\Gamma^{i}_{k})\bigg\}.
\end{flalign*}

We note that $\LPR^i$  and $\LPRS^i$ are not convex sets; however, $\LPR^i$ ($\LPRS^i$) contains to the set of extreme points of the convex set $\LCO^i$ ($\LCOS^i$). This fact plays a pivotal role later on in our analysis.

In the following, we first present two lemmas on convexity, compactness, and relationships between the preceding set of exchangeable randomized policies. We use the following lemmas, for our main results in Theorems \ref{the:1}--\ref{the:6}. 

\begin{lemma}\label{the:0}
Consider the sets of randomized policies $\LEXN^{i}, \LCOSN^{i}, \LEX^{i}$, and $\LCOS^{i}$. Then, the following holds:
\begin{itemize}
    \item [(i)] Any infinitely-exchangeable randomized policy $P_{\pi}$ in $\LEX^{i}$ also belongs to the set of randomized policies $\LCOS^{i}$, i.e., $\LCOS^{i}=\LEX^{i}$.
    \item [(ii)] $\LCOSN^{i} \subseteq \LEXN^{i}$; however, in general $\LCOSN^{i} \not= \LEXN^{i}$.
\end{itemize}
\end{lemma}

 \begin{proof}
 Part (i) has been established in \cite[Theorem 1]{SSYdefinetti2020} using a de Finetti representation theorem \cite[Theorem 1.1]{kallenberg2006probabilistic}. Part (ii) follows from the fact that finite exchangeable random variables are not necessarily conditionally i.i.d.; see e.g., \cite[Example 1.18]{schervish2012theory}, and \cite{diaconis1980finite,aldous2006ecole}.
\end{proof}

The proof of the following lemma is provided in Appendix \ref{app:0}.
\begin{lemma}\label{the:00}
Let $\mathbb{U}^{i}$ be compact, then $\LEXN^{i}$ and $\LEX^{i}$ are convex and compact.
\end{lemma}


\subsection{An existence result for a Nash equilibrium for \PN: the role of common randomness.}

In this section, we establish an existence result for general stochastic game \PN\ without imposing exchangeability on the cost function (Assumptions \ref{assump:exccost}). As in \PN\ when team plays against another team, existence of a NE may
not hold in general. This is due to lack of convexity of the space of policies for teams taking part in the game (see \cite{anantharam2007common} for an
example where such a situation arises). Considering the common randomness among DMs within teams is crucial for establishing the existence of a NE for these games since allowing the common randomness lead to the convexity of the policy space for the teams.

First, the expected costs for \PN\ under a randomized policy profile  $(P_{\pi}^{1},P_{\pi}^{2}) \in L_{N}^{1}\times  L_{N}^{2}$ are given by
\begin{flalign}
J_{\pi, N}^{i}(P_{\pi}^{1},P_{\pi}^{2})&:=\int P_{\pi}^{1}(d\underline{\gamma}^{1})P_{\pi}^{2}(d\underline{\gamma}^{2})\mu^{N}(d\omega_{0},d\underline{y})c_{N}^{i}(\underline{\gamma},\underline{y}, \omega_{0})\label{eq:finitecost}\\
&:=\int  \int c^{i}(\omega_{0},u^{1}_{1:N_{i}},u^{2}_{1:N_{i}})\prod_{i=}^{2}\prod_{k=1}^{N_{i}}\gamma_{k}^{i}(du^i_{k}|y^i_{k}) \nonumber\\
&\qquad \times  \prod_{i=1}^{2}P_{\pi}^{i}(d\gamma^{i}_{1:N_{i}})\mu^{N}(d\omega_{0}, dy^{1}_{1:N_{1}},dy^{2}_{1:N_{2}})\nonumber,
\end{flalign}
where $\mu^{N}$ is the joint probability measure of $\omega_{0}, dy^{1}_{1:N_{1}},dy^{2}_{1:N_{2}}$ and $$c_{N}^{i}(\underline{\gamma},\underline{y}, \omega_{0}):=\int c^{i}(\omega_{0},u^{1}_{1:N_{i}},u^{2}_{1:N_{i}})\prod_{i=1}^{2}\prod_{k=1}^{N_{i}}\gamma_{k}^{i}(du^i_{k}|y^i_{k}).$$


In the theorem, we establish existence of a NE for \PN. The proof of the theorem is provided in Appendix \ref{app:general}.

\begin{theorem}\label{the:general}
Consider the game \PN\ with a given IS. Let Assumptions \ref{assump:ind} hold. Suppose further that $\mathbb{U}^{i}$ is compact and $c^{i}(\omega_{0}, \cdot, \ldots, \cdot)$ is continuous and (uniformly) bounded for all $\omega_{0}$. Then,
there exists a NE $(P_{\pi}^{1*}, P_{\pi}^{2*})\in \LCON^1 \times \LCON^{2}$ for \PN\ among all policies in $L^{1}_{N} \times L^{2}_{N}$.
\end{theorem}

In Theorem \ref{the:general}, considering the common and independent randomness in policies, lead to the convexity of $\LCON^{i}$. Since $\LCON^{i}$ and $L^i_{N}$ are identical, $\LCON^{i}$ is closed. As a result, continuity of the costs together with the compactness of the actions spaces allows us to utilize Kakutani-Fan-Glicksberg fixed point theorem \cite[Corollary 17.55]{InfiniteDimensionalAnalysis}.  A similar analysis has been used for zero-sum games among two teams in \cite{hogeboom2022zero} to establish existence of a randomized Saddle-point equilibrium with common and independent randomness. Allowing a common randomness in policies may be too restrictive due to the arbitrary nature of distributions and since randomness must be externally provided.  As the number of DMs drives to infinity for \PIN, we relax this assumption, which further motivates the study of exchangeable games among teams. Additionally, in \cite{hogeboom2022zero}, it has been shown that the absolute continuity in Assumption \ref{assump:ind} among all DMs can be relaxed for zero-sum games by instead imposing Assumption \ref{assump:ind} only on the marginal distributions of observations among DMs within each team, fixing the policies of other teams. Hence, this is also true for zero-sum version of \PN\ studied here. 

\subsection{Existence of an exchangeable Nash Equilibrium for exchangeable stochastic games \PN.}

In the following, we show that the preceding class of exchangeable stochastic games admits a (finite) exchangeable NE, belonging to $ \LEXN^{1} \times  \LEXN^{2}$ (this NE is not necessarily symmetric). We first introduce the following assumptions on DMs' actions spaces, costs, and observations.

\begin{assumption}\label{assump:obsx}
For all $i=1,2$,
\begin{itemize}
    \item [(i)] $\mathbb{U}^{i}$ is compact.
    \item [(ii)] $c^{i}(\omega_{0}, \cdot, \ldots, \cdot)$ is continuous and (uniformly) bounded for all $\omega_{0}$.
    \item [(iii)] $(y^{i}_{k})_{k\in \N_{i}}$ are $N_{i}$-exchangeable, conditioned on $\omega_{0}$.
    \item [(iv)] $(y^{1}_{k})_{k\in \N_{1}}$ and $(y^{2}_{k})_{k\in \N_{2}}$ are mutually independent, conditioned on $\omega_{0}$.
    \end{itemize}
\end{assumption}

We note that since i.i.d. random variables are (infinitely) exchangeable, Assumption \ref{assump:oc} implies Assumptions \ref{assump:obsx} (iii) and (iv). The following theorem establishes existence of an exchangeable NE of \PN. The proof of the theorem is provided in Appendix \ref{app:A}.

\begin{theorem}\label{the:1}
Consider the game \PN\ with a given IS. Let Assumptions \ref{assump:exccost},  \ref{assump:ind} and \ref{assump:obsx} hold. Then,
there exists an exchangeable NE $(P_{\pi}^{1*}, P_{\pi}^{2*})$ for \PN\ among all policies in $L^{1}_{N} \times L^{2}_{N}$, i.e., $(P_{\pi}^{1*}, P_{\pi}^{2*})\in \LEXN^{1} \times  \LEXN^{2}$.
\end{theorem}

In Theorem \ref{the:1}, by continuity of the costs in Assumption \ref{assump:obsx} and convexity and compactness of $\LEXN^{i}$, we first invoke the Kakutani-Fan-Glicksberg fixed point theorem \cite[Corollary 17.55]{InfiniteDimensionalAnalysis} to show the existence of an NE among all exchangeable policy profiles. Then, we use exchangeability of the costs and observations in Assumption  \ref{assump:obsx} to show that we can restrict the search over exchangeable policies in finding best response policies to fixed exchangeable policy of other player.

Two remarks are in order: First, Theorem \ref{the:1} only guarantees the existence of an exchangeable NE which might not be symmetric since not all finite exchangeable random variables are i.i.d. (see Lemma \ref{the:0} (ii)). We cannot guarantee the existence of a symmetric NE since  restricting to symmetric  policies in finding best response policies might be with a loss for the DMs within the player. This is because, fixing policies of other players to symmetric or exchangeable policies, DMs within a deviating player faces a team problem that might not admit a symmetric optimal solution; e.g., see \cite[Example 1]{SSYdefinetti2020}. Second, in contrast to \cite{hogeboom2022zero}, in Theorem \ref{the:1}, our focus is on exchangeable games among teams, and we establish exchangeability of a NE which still might require a common randomness among DMs with teams.

\subsection{Existence of a symmetric Nash equilibrium for \PIN.}

The expected costs for \PN\ with the mean-field interaction under a randomized policy profile $(P_{\pi}^{1},P_{\pi}^{2}) \in L_{N}^{1}\times  L_{N}^{2}$ are given by
\begin{flalign}
J_{\pi, N}^{i}(P_{\pi}^{1},P_{\pi}^{2})&\::=\int P_{\pi}^{1}(d\underline{\gamma}^{1})P_{\pi}^{2}(d\underline{\gamma}^{2})\mu^{N}(d\omega_{0},d\underline{y})c_{N}^{i}(\underline{\gamma},\underline{y}, \omega_{0})\label{eq:pf}
\end{flalign}
with 
$$c_{N}^{i}(\underline{\gamma},\underline{y}, \omega_{0}):=\int \frac{1}{N_i}\sum_{k=1}^{N_{i}}c^{i}\left(\omega_{0},u^{i}_{k},\Xi^{1}(\frac{1}{N_{1}}\sum_{p=1}^{{N}_{1}}\delta_{u^{1}_{p}}), \Xi^{2}(\frac{1}{N_{2}}\sum_{p=1}^{{N}_{2}}\delta_{u^{2}_{p}})\right)\prod_{i=1}^{2}\prod_{k=1}^{N_{i}}\gamma_{k}^{i}(du^i_{k}|y^i_{k}).$$
 We note that the above cost is a special case of \PN\ introduced in \eqref{eq:finitecost} since we have a special structure for the cost function $c_{N}^{i}$ which satisfies Assumption \ref{assump:exccost}. 

Next, the expected costs for \PIN\ under a randomized policy profile $(P_{\pi}^{1},P_{\pi}^{2})\in L^{1}\times L^{2}$ are given by
\begin{flalign}
&J_{\pi, \infty}^{i}(P_{\pi}^{1},P_{\pi}^{2}):={\limsup\limits_{N_{1},N_{2} \to \infty}\int P_{\pi, N_1}^{1}(d\underline{\gamma}^{1})P_{\pi, N_2}^{2}(d\underline{\gamma}^{2})\mu^{N}(d\omega_{0},d\underline{y})c_{N}^{i}(\underline{\gamma},\underline{y}, \omega_{0})}\label{eq:pinf},
\end{flalign}
where $P_{\pi,N_i}^{i}$ is the marginal of the $P_{\pi}^{i}\in L^{i}$ to the first $N_{i}$ components, and $\mu^{N}$ is the marginal of the fixed probability measure on $(\omega_{0}, y^{1}_{1:N_{1}},y^{2}_{1:N_{2}})$.

The following theorem establishes existence of a NE for \PIN\ that is symmetric and independent. The proof of the theorem is provided in Appendix \ref{app:B}.

\begin{theorem}\label{the:2}
Consider the game \PIN\ with a given IS. Under Assumptions \ref{assump:ind}, \ref{assump:oc}, and \ref{assump:cc}, there exists an independently randomized symmetric NE profile $(P_{\pi}^{1*}, P_{\pi}^{2*})$, for \PIN\ among all policies in $L^{1} \times L^{2}$, i.e., $(P_{\pi}^{1*}, P_{\pi}^{2*}) \in \LPRS^{1} \times \LPRS^{2}$. 
\end{theorem}

We shed light into the proof of Theorem \ref{the:2}. First, since observations are conditionally i.i.d., the costs are continuous (by Assumptions \ref{assump:oc} and \ref{assump:cc}), by restricting the search to policies belonging to $\LPRS^{1}\times\LPRS^{2}$, we show that without loss we can consider a representative DM for each player (team), and then, we use the Kakutani-Fan-Glicksberg fixed point theorem to show that there exists a symmetric NE for \PIN\ among policies belonging to $\LPRS^{1}\times\LPRS^{2}$. Then, we use a de Finetti representation theorem, an argument used in \cite[Lemma 2]{SSYdefinetti2020} based on the de Finetti theorem for finite exchangeable random variables, and lower semi-continuity of the expected cost functions to show that for \PIN\, we can restrict the search in finding best response policies to fixed symmetric policy of the other player over only symmetric policies. This guarantees that an independently randomized  symmetric NE for \PIN\ constitutes a NE among all policies belonging to $L^{1}\times L^{2}$.


We note that in contrary to Theorem \ref{the:1} where we established the existence of an exchangeable randomized policies, in Theorem \ref{the:2}, we established existence of a symmetric randomized NE with only independent randomness. The reason is because, the set of randomized policies with only independent randomness in not convex subset of $L_{N}^{i}$, and hence, the Kakutani-Fan-Glicksberg fixed point theorem does not apply. However in the limit as the number of DMs drives to infinity, thanks to the \cite[Lemma 2 and Theorem 1]{SSYdefinetti2020}, we consider a single DM as a representative DM for each team with countably infinite number of DMs whose apply a symmetric policy. This allows us to use the Kakutani-Fan-Glicksberg fixed point theorem on the corresponding set of randomized policies with independent randomization (for a single DM), which is convex and compact. 

\subsection{Approximations of a symmetric Nash equilibrium for \PN}\label{eq:st-approx.}

By Theorem \ref{the:1}, there exists an exchangeable NE for \PN\ with exchangeable costs, e.g., those with the  mean-field interaction. This NE is not necessarily symmetric and varies by the number of DMs in the game. On the other hand, Theorem \ref{the:2} established the existence of a NE for \PIN\ that is symmetric and independent. In view of this, we address the following question: \emph{Does there exist a scalable (with respect to the number of DMs) symmetric approximate NE for  \PN\ with the mean-field interaction?} We answer this question in affirmative by showing that a symmetric NE for \PIN\ constitutes an approximate NE for \PN. The proof of the following theorem is provided in Appendix \ref{app:E}.

\begin{theorem}\label{the:5}
Consider the games \PN\ and \PIN\ with a given IS. Let Assumptions \ref{assump:ind}, \ref{assump:oc}, and \ref{assump:cc} hold. Then an independently randomized symmetric NE for \PIN\ constitutes an $\epsilon_{N_{1:2}}$-NE for \PN\ among all policies in $L^{1}_{N}\times L^{2}_{N}$, where $\epsilon_{N_{1:2}}\to 0$ as $N_{1}, N_{2}\to \infty$.   
\end{theorem}

Theorem \ref{the:5} is a by product of Theorem \ref{the:1} and our analysis in the proof of Theorem \ref{the:2}. We first show that a symmetric NE for \PIN\ constitutes an $\epsilon_{N_{1:2}}$-NE for \PN\ among symmetric ones. Then, as an implication of analysis developed in the proof of Theorem \ref{the:2}, using Theorem \ref{the:1} and Lemma \ref{the:0}, we establish that a symmetric NE for \PIN\ constitutes an $\epsilon_{N_{1:2}}$-NE for \PN\ among all randomized policies.

\section{Exchangeable Dynamic Games with Symmetric Information Structures.}\label{sec:4}
In this section, we study exchangeable dynamic stochastic games.  We first introduce randomized policies with their suitable topology, and then, we establish existence of a NE for \PNT\ that is exchangeable, and a NE for \PIT\ that is symmetric. 

\subsection{Topology on control policies for dynamic games.}

We first introduce two reduction conditions as variations of Assumption \ref{assump:ind} for static games, adapted to the dynamic game setting with a given IS, that allow us to define randomized policies for dynamic games  as Borel spaces.

\begin{assumption}\label{assump:ind1}
Let for every $N\in \mathbb{N}\cup\{\infty\}$, $\nu_{t}^{N}$ be the distribution of $\underline{y}_{t}^{1:2}$, conditioned on the information history $H_{t}:=\omega_{0}, \underline{x}_{0}^{1:2}, \underline\zeta_{0:t-1}^{1:2}, \underline{y}^{1:2}_{0:t-1}, \underline{u}_{0:t-1}^{1:2}$. 
One of the following conditions holds for every $N\in \mathbb{N}\cup\{\infty\}$ and for every DM at $t=0,\dots,T-1$:
\begin{itemize}
\item [(i)] (Independent reduction):  There exist a probability measure $\tau_{k,t}^{i}$ and a function $\psi_{k,t}^{i}$ such that for all Borel sets $A^{i}_{k}$ on $\mathbb{Y}^{i}$ (with $A=\prod_{i=1}^{2}A^i_{1} \times \dots \times A^i_{N_{i}}$)
\begin{flalign*}
&\nu_{t}^{N}(A|H_{t})=\prod_{i=1}^{2}\prod_{k=1}^{N_{i}}\int_{A^i_{k}}\psi_{k, t}^{i}(y^{i}_{k, t}, H_{t})\tau_{k,t}^{i}(dy_{k, t}^{i}).
\end{flalign*}
\item [(ii)] (Nested reduction): There exist a probability measure $\eta^{i}_{k,t}$ and a function $\phi^{i}_{k,t}$ such that for all Borel sets $A^{i}_{k}$ on $\mathbb{Y}^{i}$ 
\begin{flalign*}
\nu_{t}^{N}(A|H_{t})=\prod_{i=1}^{2}\prod_{k=1}^{N_{i}}\int_{A^i_k}\phi_{k,t}^{i}(y^{i}_{k,t}, H_{t})\eta_{k,t}^{i}(dy_{k,t}^{i}|H_{k,t}^{i})\nonumber,
\end{flalign*}
with $H_{k,t}^{i}:={x}_{k,0}^{i}, \zeta_{k,0:t-1}^{i}, {y}^{i}_{k,0:t-1}, {u}_{k,0:t-1}^{i}$, and for each DM$^{i}_{k}$ over time, there exists a static reduction under which the IS of DMs is expanding, i.e., $\sigma(y_{k,t}^{i}) \subset \sigma(y_{k,t+1}^{i})$.
\end{itemize}
\end{assumption}

Let $\mathbb{T}:=0,\ldots,T-1$. Using Assumption \ref{assump:ind1} (i), we equip randomized policies with a suitable topology. This topology yields that 
\begin{flalign*}
\pmb{\gamma^{i, n}_{k}}  \xrightarrow[n \to \infty]{} \pmb{\gamma^{i}_{k}} \iff {\gamma}^{i,n}_{k,t}(d{u}^{i}_{k,t}|{y}^{i}_{k,t})\tau^{i}_{k,t}(d{y}^{i}_{k,t}) \xrightarrow[n \to \infty]{\text{weakly}} {\gamma}^{i}_{k,t}(d{u}^{i}_{k,t}|{y}^{i}_{k,t})\tau^{i}_{k,t}(d{y}^{i}_{k,t})\quad \forall~t\in \mathbb{T}.
\end{flalign*}
Similarly, under Assumption \ref{assump:ind1}(ii), the topology on randomized policies yields that
 \begin{flalign*}
\pmb{\gamma^{i, n}_{k}}  \xrightarrow[n \to \infty]{} \pmb{\gamma^{i}_{k}} \iff {\gamma}^{i, n}_{k,t}(d{u}^{i}_{k,t}|{y}^{i}_{k,0:t})\eta^{i}_{k,t}(d{y}^{i}_{k,0:t}) \xrightarrow[n \to \infty]{\text{weakly}} {\gamma}^{i}_{k,t}(d{u}^{i}_{k,t}|{y}^{i}_{k,0:t})\eta^{i}_{k,t}(d{y}^{i}_{k,0:t})\quad \forall~t\in \mathbb{T}.
\end{flalign*}
 That is, under Assumption \ref{assump:ind1}, by considering $\pmb{\gamma^{i}_{k}}$, we are able to define randomized policies similar to those in Section \ref{sec:strategic} adapted to the dynamic setting\footnote{For examples under which reductions in Assumption \ref{assump:ind1} hold, see \cite[p. 15]{SSYdefinetti2020} and \cite{SBSYgamesstaticreduction2021}.}. 

\subsection{Existence of an exchangeable  Nash equilibrium for \PNT.}

Here, we establish existence of an exchangeable NE for a class of symmetric games general than \PNT. We study a general formulation for \PNT\ where the ISs is  $I_{k,t}^{i}:=\{y_{k,t}^{i}\}$, and is symmetric given by \begin{align}\label{eq:syminfo}
y^{i}_{k,t}=h_{t}^{i}(x_{k,0}^{i},\zeta^{i}_{k,0:t},  u_{k,0:t-1}^{i}, M_{-k, 0:t-1}^{1:2}),
\end{align}
where for all $k\in \mathcal{N}_{i}$, the functions $h_{t}^{i}$ are identical, and we used $$M_{-k, 0:t-1}^{1:2}:=(x_{-k,0}^{i},\zeta^{i}_{-k,0:t},  u_{-k,0:t-1}^{i},\underline{x}_{0}^{-i},\underline\zeta^{-i}_{0:t},  \underline{u}_{0:t-1}^{-i})$$ to denote the information up to time $t$ excluding $(x_{k,0}^{i},\zeta^{i}_{k,0:t},  u_{k,0:t-1}^{i})$. We note that the ISs of \PNT\ and  \PIT, given by \eqref{eq:mfdynamics2} and \eqref{eq:mfobs2}, are symmetric and special cases of \eqref{eq:syminfo}.

 The expected costs within  $(P_{\pi}^{1},P_{\pi}^{2})\in L_{N}^{1} \times L_{N}^{2}$ are given by
\begin{flalign}
J_{\pi, N}^{i,T}(P_{\pi}^{1},P_{\pi}^{2})&:=\int  P_{\pi}^{1}(d\pmb{\underline{\gamma}^{1}})P_{\pi}^{2}(d\pmb{\underline{\gamma}^{2}})\mu^{N}(d\omega_{0},d\pmb{\underline{\zeta}})c_{N}^{i}(\pmb{\underline{\gamma}},\pmb{\underline{y}}, \omega_{0})\nu^{N}(d\pmb{\underline{y}}|\pmb{\underline{\zeta}}, \pmb{\underline{\gamma}},\pmb{\omega_{0}})\label{eq:finitecostd}\\
&:=\int \left(\int c^{i}(\omega_{0},\pmb{\underline{\zeta}^{1:2}},\pmb{\underline{u}^{1:2}})\prod_{i=1}^{2}\prod_{k=1}^{N_{i}}\pmb{\gamma^{i}_{k}}(d\pmb{u^{i}}|\pmb{y^{i}})\right)\nonumber\\
&\qquad \times P_{\pi}(d\pmb{\underline\gamma^{1:2}})\mu^{N}(d\omega_{0}, d\pmb{\underline\zeta^{1:2}}) \prod_{t=0}^{T-1}{{\nu_{t}^{1:2}}}\left(d\underline{y}^{1:2}_{t}\middle|H_{t}\right)\nonumber,
\end{flalign}
where $\mu^{N}$ denotes the joint distribution of $\omega_{0}, \pmb{\underline{\zeta}^{1:2}}$, and
$$c_{N}^{i}(\pmb{\underline{\gamma}},\pmb{\underline{y}}, \omega_{0}):=\int c^{i}(\omega_{0},\pmb{\underline{\zeta}^{1:2}},\pmb{\underline{u}^{1:2}})\prod_{i=1}^{2}\prod_{k=1}^{N_{i}}\pmb{\gamma^{i}_{k}}(d\pmb{u^{i}}|\pmb{y^{i}}),$$ satisfies the following exchangeability assumption.

\begin{assumption}\label{assump:3.1}
$c^{i}$ is (separately) exchangeable  for any permutations $\sigma$ and $\tau$ of $\{1,\ldots,N_{1}\}$ and $\{1,\ldots,N_{2}\}$, respectively, for all $\omega_{0}, \pmb{\underline{\zeta}^{1:2}}$ \begin{flalign}\label{eq:3.3}
c\left(\omega_{0},\pmb{\zeta_{\sigma(1):\sigma(N_i)}^{1}},\pmb{\zeta_{\tau(1):\tau(N_{i})}^{2}}, \pmb{u_{\sigma(1):\sigma(N_{i})}^{1}},\pmb{u_{\tau(1):\tau(N_{i})}^{2}}\right)=c^{i}\left(\omega_{0}, \pmb{\underline{\zeta}^{1:2}},\pmb{\underline{u}^{1:2}}\right).
 \end{flalign}
\end{assumption}

In addition, we let the following assumptions hold. 
\begin{assumption}\label{assump:sxcdy}
For $i=1,2$, 
\begin{itemize}
\item [(i)] $\mathbb{U}^{i}$ is compact.
\item [(ii)]  $c^{i}(\omega_{0}, \pmb{\underline{\zeta}^{1:2}}, \cdot)$ in \eqref{eq:finitecostd}, is continuous and (uniformly) bounded for all $\omega_{0}, \pmb{\underline{\zeta}^{1:2}}$.
\item [(iii)] $h_{t}^{i}$ in \eqref{eq:syminfo}, is continuous in actions for $t=0, \ldots, T-1$.
\item [(iv)] $\pmb{\underline\zeta^{1:2}}$ are exchangeable, conditioned on $\omega_{0}$; 
\item  [(v)] For all Borel sets $A^{i}_{k}$ on ${{\mathbb{Y}^{i}}}$ (with $A^{i}=A^{i}_{1}\times \dots \times A^{i}_{N_{i}}$) 
\begin{flalign}
&{\nu_{t}^{N}}\left(A\middle|H_{t}\right)=\prod_{i=1}^{2}\prod_{k=1}^{N_{i}}{\nu_{t}^{i}}\left(A^{i}\middle|\omega_{0}, x_{k,0}^{i},{\zeta}^{i}_{k,0:t-1},\underline{y}_{0:t-1}^{1:2},\underline{u}^{1:2}_{0:t-1}\right),\quad \forall t\in \mathbb{T}\nonumber.
\end{flalign}
\end{itemize}

\end{assumption}

The costs of \PNT \ and \PIT, given in \eqref{eq:mfcost}, satisfy Assumption \ref{assump:3.1}. Also, for the state dynamics and observations of  \PNT \ and \PIT, Assumption \ref{assump:2} implies Assumption \ref{assump:sxcdy}. In the following theorem, we establish the existence of an exchangeable NE for exchangeable dynamic games \PNT\ with a symmetric IS \eqref{eq:syminfo}. The proof of the theorem is provided in Appendix \ref{app:C}. 

\begin{theorem}\label{the:3}
Consider the game \PNT \ with a given symmetric IS. Let Assumptions \ref{assump:ind1},  \ref{assump:3.1} and \ref{assump:sxcdy} hold. Then,
there exists an exchangeable NE $(P_{\pi}^{1*}, P_{\pi}^{2*})$ for \PNT\ among all policies in $L^{1}_{N} \times L^{2}_{N}$, i.e., $(P_{\pi}^{1*}, P_{\pi}^{2*})\in \LEXN^{1} \times  \LEXN^{2}$;
\end{theorem}

Similar to Theorem \ref{the:1}, invoking  the Kakutani-Fan-Glicksberg fixed point theorem \cite[Corollary 17.55]{InfiniteDimensionalAnalysis}  establishes the existence of an exchangeable NE among all exchangeable policy profiles. Since the IS is symmetric, by Assumption \ref{assump:3.1} and Assumptions \ref{assump:sxcdy} (iv) and (v), we then show that restriction to exchangeable policies is without loss for finding the best response policies to exchangeable policies of others.

\subsubsection{Existence of a symmetric Nash equilibrium for \PIT.}

In this subsection, we establish existence of a symmetric NE for \PIT.  Define state dynamics and observations as \eqref{eq:mfdynamics2} and \eqref{eq:mfobs2}, respectively. Let $I_{k,t}^{i}=\{{y}^{i}_{k,t}\}$, and $\underline\zeta_{t}^{i}:=(\underline{w}_{t}^{i}, \underline{v}_{t}^{i})$ with $\underline\zeta_{0}^{i}:=(\underline{x}_{0}^{i}, \underline{w}_{0}^{i}, \underline{v}_{0}^{i})$.

 Let the expected costs for \PNT\ within a randomized policy profile $(P_{\pi}^{1},P_{\pi}^{2})\in L_{N}^{1} \times L_{N}^{2}$ be given by
\begin{flalign}
J_{\pi, N}^{i,T}(P_{\pi}^{1},P_{\pi}^{2})&=\int  P_{\pi}^{1}(d\pmb{\underline{\gamma}^{1}})P_{\pi}^{2}(d\pmb{\underline{\gamma}^{2}})\mu^{N}(d\omega_{0},d\pmb{\underline{\zeta}})c_{N}^{i}(\pmb{\underline{\gamma}},\pmb{\underline{y}}, \omega_{0})\nu^{N}(d\pmb{\underline{y}}|\pmb{\underline{\zeta}}, \pmb{\underline{\gamma}},\pmb{\omega_{0}}),\label{eq:finiterandom}
\end{flalign}
where 
\begin{flalign*}
c^{N}(\underline{\zeta}, \underline{\gamma},\underline{y}, \omega_{0})
&:=\int \frac{1}{N_{i}}\sum_{k=1}^{N_{i}}\sum_{t=0}^{T-1}c^{i}\bigg(\omega_{0},x_{k, t}^{i},u_{k, t}^{i},\Xi^{1}_{x}(\frac{1}{N_{1}}\sum_{p=1}^{N_{1}}\delta_{x_{p, t}^{1}}), \Xi^{2}_{x}(\frac{1}{N_{2}}\sum_{p=1}^{N_{2}}\delta_{x_{p, t}^{2}}),\\
& \qquad \qquad \Xi^{1}_{u}(\frac{1}{N_{1}}\sum_{p=1}^{N_{1}}\delta_{u_{p, t}^{1}}), \Xi^{2}_{u}(\frac{1}{N_{2}}\sum_{p=1}^{N_{2}}\delta_{u_{p, t}^{2}})\bigg)\prod_{i=1}^{2}\prod_{k=1}^{N_{i}}\pmb{\gamma^{i}_{k}}(d\pmb{u^{i}}|\pmb{y^{i}}).
\end{flalign*}

{

 Let the expected costs for \PIT\ within  $(P_{\pi}^{1},P_{\pi}^{2})\in L^{1} \times L^{2}$ be given by
\begin{flalign}
&J_{\pi, \infty}^{i,T}(P_{\pi}^{1},P_{\pi}^{2})=\limsup\limits_{N \to \infty}\int  P_{\pi,N}^{1}(d\pmb{\underline{\gamma}^{1}})P_{\pi,N}^{2}(d\pmb{\underline{\gamma}^{2}})\mu^{N}(d\omega_{0},d\pmb{\underline{\zeta}})c_{N}^{i}(\pmb{\underline{\gamma}},\pmb{\underline{y}}, \omega_{0})\nu^{N}(d\pmb{\underline{y}}|\pmb{\underline{\zeta}}, \pmb{\underline{\gamma}},\pmb{\omega_{0}}),\label{eq:infiniterandom}
\end{flalign}
where $P_{\pi,N}^{i}$ is the marginal of the $P_{\pi}^{i}\in L^{i}$ to the first $N_{i}$ components, and $\mu^{N}$ is the joint distribution of $(\omega_{0}, \pmb{{\zeta}^{1}_{1:N_{1}}}, \pmb{{\zeta}^{2}_{1:N_{2}}})$.

\begin{assumption}\label{assump:ind2}
Assumption \ref{assump:ind1} holds with functions $\psi_{t}^{i}$ and $\phi_{t}^{i}$ are of the following forms for every $i\in \cal{N}$ and $t=0,\dots, T-1$: 
\begin{flalign*}
&\psi_{k,t}^i\bigg(y^{i}_{k,t}, \omega_{0}, \zeta_{k,0:t-1}^{i}, {y}^{i}_{k,0:t-1}, {u}_{k,0:t-1}^{i}, \frac{1}{N_{1}}\sum_{p=1}^{N_{1}}\delta_{u_{p,0:t-1}^{1}}, \frac{1}{N_{2}}\sum_{p=1}^{N_{2}}\delta_{u_{p,0:t-1}^{2}},\\
&\qquad\qquad\frac{1}{N_{1}}\sum_{p=1}^{N_{1}}\delta_{x_{p,0:t-1}^{1}}, \frac{1}{N_{2}}\sum_{p=1}^{N_{2}}\delta_{x_{p,0:t-1}^{2}}\bigg),\\
&\phi_{k,t}^i\left(y^{i}_{k,t}, \omega_{0}, \frac{1}{N_{i}}\sum_{p=1}^{N_{1}}\delta_{u_{p,0:t-1}^{1}}, \frac{1}{N_{2}}\sum_{p=1}^{N_{2}}\delta_{u_{p,0:t-1}^{2}}, \frac{1}{N_{1}}\sum_{p=1}^{N_{1}}\delta_{x_{p,0:t-1}^{1}}, \frac{1}{N_{2}}\sum_{p=1}^{N_{2}}\delta_{x_{p,0:t-1}^{2}}\right),
\end{flalign*}
where $\psi_{k,t}^{i}(y^{i}_{k,t}, \omega_{0}, \zeta_{k,0:t-1}^{i}, {y}^{i}_{k,0:t-1}, \cdot, \ldots, \cdot)$ and $\phi_{t}^{i}(y^{i}_{k,t}, \omega_{0}, \cdot, \ldots, \cdot)$ are continuous (weakly continuous in their four last arguments) and uniformly bounded for all $y^{i}_{k,t}$, $\omega_{0}$, $\zeta_{k,0:t-1}^{i},{y}^{i}_{k,0:t-1}$.
\end{assumption}

In the following theorem, we establish the existence of a symmetric randomized NE for \PIT. The proof of the theorem is provided in Appendix \ref{app:D}. 

\begin{theorem}\label{the:4}
Consider the game \PIT\ with a given IS. Let Assumptions \ref{assump:c}, \ref{assump:2}, and \ref{assump:ind2} hold. Then, there exists an independently randomized symmetric Nash  equilibrium profile $(P_{\pi}^{1*}, P_{\pi}^{2*})$ for \PIT\ among all policies in $L^{1} \times L^{2}$, i.e., $(P_{\pi}^{1*}, P_{\pi}^{2*})\in \LPRS^1 \times \LPRS^2$.
\end{theorem}

Theorem \ref{the:4} follows from an argument similar to that used for Theorem \ref{the:2} with additional technical argument using continuity and exchangeability of the cost, state dynamics and observations.

\subsection{Approximations of a symmetric Nash equilibrium for \PNT.}
Similar to the static setting, we show that there exists a symmetric approximate NE for \PNT with the mean-field interaction, by showing that a symmetric NE for \PIN\ constitutes an approximate NE for \PN. 

\begin{theorem}\label{the:6}
Consider the games \PNT\ and \PIT\ with a given IS. Let Assumptions \ref{assump:c}, \ref{assump:2}, and \ref{assump:ind2} hold. Then any independently randomized symmetric NE for \PIT\ constitutes an $\epsilon_{N_{1:2}}$-NE for \PNT\ among all policies in $L^{1}_{N}\times L^{2}_{N}$, where $\epsilon_{N_{1:2}}\to 0$ as $N_{1}, N_{2}\to \infty$.  
\end{theorem}

\begin{proof}
The proof follows from an argument similar to the proof of Theorem \ref{the:5}.
\end{proof} 

\section{Conclusion.}
In this paper, we have studied stochastic static and dynamic mean-field games among teams with finite and infinite number of decision makers. We have established existence of a NE for exchangeable games among teams with finite number of DMs, and have shown that this NE is exchangeable. For mean-field games among teams with infinite number of decision makers, we have established existence of a NE, and have shown that it is symmetric. Finally, we have established existence of an approximate NE that exhibits symmetry for games among teams with finite number of decision makers and scalable (with respect to the number of decision makers of each team), using a symmetric NE of the corresponding mean-field game.

\section{Appendix.}

\subsection{Proof of Lemma \ref{the:00}.}\label{app:0}
 The set $\LEXN^{i}$ is a non-empty convex subset of locally convex set $L_{N}^{i}$. Since $\mathbb{U}^{i}$ is compact, the marginal of probability measures on $\mathbb{U}^{i}$ is tight. Since the probability measure on $\mathbb{Y}^{i}$ is fixed, the marginal on $\mathbb{Y}^{i}$ is also tight. Since marginals are tight, the collection of all measures on $\mathbb{U}^{i} \times \mathbb{Y}^{i}$ with these tight marginals is also tight (see e.g., \cite[Proof of Theorem 2.4]{yukselSICON2017}), and hence, the set $\Gamma^{i}_{k}$ is tight for each $i=1,2$ and $k \in \N_{i}$. Hence, $\LEXN^{i}$ is tight for $i=1,2$. Next, we show that $L_{\text{EX},N}^{i}$ is closed under the weak-convergence topology. Suppose that $P_{\pi}^{i, \infty}$ is the limit, in the weak convergence topology, of a converging sequence of randomized policies $\{P_{\pi}^{i, n}\}_{n}$ as $n \to \infty$. Also, suppose that $P_{\pi}^{\sigma, i, \infty}$ is the limit, in the weak convergence topology, of a converging sequence of randomized policies $\{P_{\pi}^{\sigma, i, n}\}_{n}$ as $n \to \infty$, where for $A^{i}_{k} \in \mathcal{B}(\Gamma^{i}_{k})$ and for  all $\sigma \in S_{N_{i}}$
\begin{flalign*}
&P_{\pi}^{\sigma, i, n}(\gamma^{i}_{1}\in A^{i}_{1},\cdots, \gamma^{i}_{N_{i}}\in A^{i}_{N_{i}})=P_{\pi}^{i, n}(\gamma^{i}_{\sigma(1)}\in A^{i}_{1},\cdots, \gamma^{i}_{\sigma(N_{i})}\in A^{i}_{N_{i}}).
\end{flalign*}

 Let $\mathcal{T}$ be a countable measure-determining subset of the set of all real-valued continuous functions on $\prod_{k=1}^{N_i}\Gamma^{i}_k$. For a function $f\in \mathcal{T}$, we have
\begin{align}
    \int f(\gamma^{i}_{1:N_{i}})dP_{\pi}^{\sigma,i, \infty}
    &= \lim_{n\to \infty}\int f(\gamma^{i}_{1:N_{i}})dP_{\pi}^{i, n}
    =\int f(\gamma^{i}_{1:N_{i}})dP_{\pi}^{i, \infty}\label{eq:30-the1},
\end{align}
where \eqref{eq:30-the1} follows from the fact that $P_{\pi}^{i, n}$ is $N_{i}$-exchangeability for every $n$, and the weak convergence of   $\{P_{\pi}^{\sigma,i, n}\}_{n}$ and $\{P_{\pi}^{i, n}\}_{n}$. Since $\mathcal{T}$ is countable, for all $f\in \mathcal{T}$, we get \eqref{eq:30-the1}.
Since $\mathcal{T}$ is measure-determining,  we get that $P_{\pi}^{\sigma,i, \infty}$ and $P_{\pi}^{i, \infty}$ are the same, and hence, $\LEXN^{i}$ is compact as it is tight. Convexity of $\LEXN^{i}$ follows from the fact that the convex combination of exchangeable probability measures remains an exchangeable probability measure. An argument similar to the above establishes that $\LEX^{i}$ is convex and compact.


\subsection{Proof of Theorem \ref{the:general}.}\label{app:general}
We use the Kakutani-Fan-Glicksberg fixed point theorem \cite[Corollary 17.55]{InfiniteDimensionalAnalysis}. By \cite[Theorem A.1]{SSYdefinetti2020}, the set $\LCON^{i}$ and $L_{N}^{i}$ are identical. This is due to the fact that $\LCON^{i}$ is a convex set and the set of extreme points of $L_{N}^{i}$ is a subset of $\LCON^{i}$. This implies that the set $\LCON^{i}$ is a non-empty, convex, and closed. Using a similar argument as that in the proof of Lemma \ref{the:00}, since the action spaces are compact, $\LCON^{i}$ is compact since it is tight. 

 Define the best response correspondence $\Phi:\prod_{i=1}^{2}\LCON^{i} \to 2^{\prod_{i=1}^{2}\LCON^{i}}$ as $\Phi(P_{\pi}^{1}, P_{\pi}^{2}):=\prod_{i=1}^{2} \text{BR}(P_{\pi}^{-i})$, where for  $i=1, 2$,
\begin{flalign}
&\text{BR}(P_{\pi}^{-i}) :=\left\{P_{\pi}^{i}\in \LCON^{i} \bigg| J^{i}_{\pi}({P}_{\pi}^{i}, P_{\pi}^{-i}) \leq J^{i}_{\pi}(\hat{P}_{\pi}^{i}, P_{\pi}^{-i})  \:\:\forall \hat{P}_{\pi}^{i}\in  \LCON^{i}\right\}\label{ea:BR1}.
\end{flalign}
In the following, we show that $\Phi$ admits a fixed point, which will be denoted by $(P_{\pi}^{1*}, P_{\pi}^{2*})$.  

Next, we show that the graph $$G=\left\{\left((P_{\pi}^{1}, P_{\pi}^{2}),  \Phi(P_{\pi}^{1}, P_{\pi}^{2})\right)\bigg| (P_{\pi}^{1}, P_{\pi}^{2}) \in \prod_{i=1}^{2}\LCON^{i}\right\}$$ is closed. If a sequence of policies $\{{P}_{\pi}^{-i, n}\}_{n}\subseteq \LCON^{-i}$  converges weakly to ${P}_{\pi}^{-i}$, and if ${P}_{\pi}^{i, n}\in \text{BR}({P}_{\pi}^{-i, n})$, we get ${P}_{\pi}^{i}\in \text{BR}({P}_{\pi}^{-i})$ since
\begin{flalign}
\inf_{\hat{P}_{\pi}^{i}\in \LCON^{i}}J_{\pi}^{i}(\hat{P}_{\pi}^{i}, {P}_{\pi}^{-i}) &= \inf_{\hat{P}_{\pi}^{i}\in \LCON^{i}}\lim_{n\to \infty}J_{\pi}^{i}(\hat{P}_{\pi}^{i}, {P}_{\pi}^{-i, n}) \label{eq:14}\\
&\geq J_{\pi}^{i}(P_{\pi}^{i}, {P}_{\pi}^{-i})\label{eq:16}\\
&\geq \inf_{\hat{P}_{\pi}^{i}\in \LCON^{i}}J_{\pi}^{i}(\hat{P}_{\pi}^{i}, {P}_{\pi}^{-i})\label{eq:17}.
\end{flalign}
Equality \eqref{eq:14} follows from the generalized dominated convergence theorem for varying measures \cite[Theorem 3.5]{serfozo1982convergence} since the bounded costs $c^{i}$ converge continuously\footnote{We recall that the sequence $\{f_{n}\}_{n}$ converges continuously to $f$ if and only if $f_{n}(a_{n})\to f(a)$ whenever $a_{n}\to a$ as $n \to \infty$.}, and  $\{{P}_{\pi}^{-i, n}\}_{n}$ converges weakly to ${P}_{\pi}^{-i}$. Inequality \eqref{eq:16} follows from exchanging limit and infimum, ${P}_{\pi}^{i, n}\in \text{BR}({P}_{\pi}^{-i, n})$, and an argument similar to that used for \eqref{eq:14}, and \eqref{eq:17} follows from the fact that $P_{\pi}^{i}\in \LCON^{i}$ (since $\LCON^{i}$ is closed). This implies that the graph $G$ is closed. Moreover, since $\LCON^{i}$ is compact, and $c$ are continuous in actions, using \cite[Proposition D.5(b)]{HernandezLermaMCP}, the map $$F^{i}:P_{\pi}^{-i} \mapsto \inf_{\hat{P}_{\pi}^{i}\in \LCON^{i}}J_{\pi}^{i}(\hat{P}_{\pi}^{i}, {P}_{\pi}^{-i})$$ is
continuous. This implies that $\text{BR}({P}_{\pi}^{-i})$ is non-empty. 
 Thus, by the Kakutani–Fan–Glicksberg fixed point theorem, $\Phi$ admits a fixed point, and this completes the proof since the set $\LCON^{i}$ and $L_{N}^{i}$ are identical.

\subsection{Proof of Theorem \ref{the:1}.}\label{app:A}
By Lemma \ref{the:00}, the set $\LEXN^{i}$ is a non-empty convex compact subset of locally convex set $L_{N}^{i}$. Hence, using  Kakutani-Fan-Glicksberg fixed point theorem \cite[Corollary 17.55]{InfiniteDimensionalAnalysis} similar to the proof of Theorem \ref{the:general}, we show that there exists a NE $(P_{\pi}^{1*}, P_{\pi}^{2*})\in \LEXN^{1}\times \LEXN^{2}$ among all policies in $ \LEXN^{1} \times \LEXN^{2}$. We only need to show that $(P_{\pi}^{1*}, P_{\pi}^{2*})$ constitutes a NE  among all randomized policies in $L^{1}_{N} \times L^{2}_{N}$. It is sufficient to show that for any $P_{\pi}^{-i*}\in L_{\text{EX},N}^{-i}$, we have
    \begin{flalign}
    \inf_{\hat{P}_{\pi}^{i}\in L_{N}^{i}}J_{\pi}^{i}(\hat{P}_{\pi}^{i}, {P}_{\pi}^{-i*})=\inf_{\hat{P}_{\pi}^{i}\in \LEXN^{i}}J_{\pi}^{i}(\hat{P}_{\pi}^{i}, {P}_{\pi}^{-i*}).\label{eq:e-cost}
    \end{flalign}
    
    This follows from \cite[Lemma 1]{SSYdefinetti2020} using Assumptions \ref{assump:exccost} and \ref{assump:obsx} (iii) and (iv), by first showing that for any policy $\hat{P}_{\pi}^{i}\in L_{N}^{i}$, we have that 
    $J_{\pi}^{i}(\hat{P}_{\pi}^{i}, {P}_{\pi}^{-i*})=J_{\pi}^{i}(\hat{P}_{\pi}^{\sigma, i}, {P}_{\pi}^{-i*})$, where $\hat{P}_{\pi}^{\sigma, i}$ is the $\sigma$-permutation of $\hat{P}_{\pi}^{i}$, and then by showing that an exchangeable policy constructed as an average of all possible permutations of ${P}_{\pi}^{i}$ does not perform worst for PL$^{i}$. This implies \eqref{eq:e-cost}, and hence, under Assumption \ref{assump:ind}, there exists a NE profile $({P}_{\pi}^{1*}, {P}_{\pi}^{2*})$, belonging to $\LEXN^{1}\times \LEXN^{2}$. We have that every NE policy profile for the game under the reduction of Assumption \ref{assump:ind} constitutes a NE policy profile for the original game \cite[Theorem 3.1]{SBSYgamesstaticreduction2021} since policies do not change under the reduction  and the expected costs are the same, and ISs are equivalent. Hence, $({P}_{\pi}^{1*}, {P}_{\pi}^{2*})$ constitutes NE for the \PN, and this completes the proof.

\subsection{Proof of Theorem \ref{the:2}.}\label{app:B}
The proof proceeds in two steps. In step 1, we show that there exists a policy profile $({P}_{\pi, \infty}^{1*}, {P}_{\pi, \infty}^{2*})$ that constitutes a NE for \PIN, among all randomized policies in $\LPRS^{1} \times \LPRS^{2}$. In step 2, we show that the policy profile $({P}_{\pi, \infty}^{1*}, {P}_{\pi, \infty}^{2*})$ constitutes a NE for \PIN, among all randomized policies in $L^{1} \times L^{2}$.

\begin{itemize}[wide]
    \item [{\it Step 1.}] We first restrict the search in finding a NE for \PIN\ to randomized policies belonging to $\LPRS^{1} \times \LPRS^{2}$.  We note that under Assumptions \ref{assump:ind} and \ref{assump:oc} (see \eqref{eq:abscon-iden}), observations of DMs within players are i.i.d. and also independent of $\omega_{0}$ via using a change of measure argument in \eqref{eq:abscon-iden}. Since our search is restricted to $\LPRS^{1} \times \LPRS^{2}$,  by the strong law of large numbers and the generalized dominated convergence theorem for varying measures \cite[Theorem 3.5]{serfozo1982convergence} since the bounded costs $c^{i}$ converge continuously, we have for any $(P_{\pi}^{1},P_{\pi}^{2})$ in $\LPRS^{1} \times \LPRS^{2}$ 
    \begin{flalign}
&\limsup_{N_{1},N_{2}\to \infty}\int \frac{1}{N_{i}}\sum_{k=1}^{N_{i}}c^{i}\left(\omega_{0},u^{i}_{k}, \Xi^{1}\left(\frac{1}{N_{1}}\sum_{k=1}^{N_{1}}u^{1}_{k}\right), \Xi^{2}\left( \frac{1}{N_{2}}\sum_{k=1}^{N_{2}}u^{2}_{k}\right)\right) \nonumber\\
    &\qquad \times {P}_{\pi}^{i}(d\underline{\gamma}^{i})P_{\pi}^{-i}(d\underline{\gamma}^{-i})\mathbb{P}^{0}(d\omega_{0})\prod_{i=1}^{2}\prod_{k=1}^{N_{i}}\hat{f}^{i}(d{y}^{i}_{k}, \omega_{0})Q^{i}(d{y}^{i}_{k})\gamma_{k}^{i}(du^i_{k}|y^i_{k})\nonumber\\
&=\int c^{i}\left(\omega_{0},u^{i}_{\R},\Xi^{1}\left(\Lambda^{1}(du^{1}_{\R})\right), \Xi^{2}\left(\Lambda^{2}(du^{2}_{\R})\right)\right)\Lambda^{i}(du^{i}_{\R}) \nonumber\\
&\qquad \times \mathbb{P}^{0}(d\omega_{0})\prod_{s=1}^{2}\hat{f}^{s}\left({y}^{s}_{\R}, \omega_{0}\right)Q^{s}(dy^{s}_{\R})P_{\pi,\R}^{s}(d{\gamma}^{s}_{\R})\gamma_{\R}^{s}(du^s_{\R}|y^s_{\R})\nonumber\\
&:=\bar{J}_{\pi, \infty}^{i}(P_{\pi, \R}^{i},P_{\pi, \R}^{-i}, \Lambda^{i}, \Lambda^{-i})\nonumber,
\end{flalign}
where the sub-index \R\ denotes the representative DM for each player $i$. In the above, 
\begin{align}\label{eq:policyrep}
    P_{\pi}^{i}(\underline\gamma^{i} \in \cdot)=\prod_{k\in \mathbb{N}}{P}_{\pi, \R}^{i}(\gamma^{i}_{k} \in \cdot),
\end{align}
for some ${P}_{\pi, \R}^{i}(\gamma^{i}_{k} \in \cdot)$ belonging to the set $\mathcal{P}(\Gamma^{i}_{\R})$, and 
\begin{align}
\Lambda^{i}(\cdot)=\mathcal{L}(u^{i}_{\R}).\label{eq:consistency}
\end{align}
 In view of the above costs, a policy profile $({P}_{\pi}^{1}, {P}_{\pi}^{2})$ of the forms \eqref{eq:policyrep} constitutes a NE for \PIN\ if ${P}_{\pi,\R}^{i}\in \text{BR}({P}_{\pi,\R}^{-i}, \Lambda^{1}, \Lambda^{2})$ for $i=1,2$ with 
\begin{flalign}
&\text{BR}({P}_{\pi,\R}^{-i}, \Lambda^{1}, \Lambda^{2}):=\bigg\{{P}_{\pi,\R}^{i}\in \mathcal{P}(\Gamma^{i}_{\R})\bigg|\label{eq:BR2}\\
&\bar{J}_{\pi, \infty}^{i}({P}_{\pi, \R}^{i},P_{\pi, R}^{-i}, \Lambda^{i}, \Lambda^{-i}) \leq \bar{J}_{\pi, \infty}^{i}(\hat{P}_{\pi, \R}^{i},P_{\pi,\R}^{-i}, \Lambda^{i}, \Lambda^{-i})\quad \forall \hat{P}_{\pi,\R}^{i} \in \mathcal{P}(\Gamma^{i}_{\R})\bigg\}\nonumber,
\end{flalign}
and the consistency condition \eqref{eq:consistency} holds. We note that the above best response map is different from that defined in \eqref{ea:BR1}, this is because, in \eqref{eq:BR2}, each representative DM policy of a player is a best response to the policy of the  representative DM of the other player and the mean-field terms of both players.

Define the best response correspondence $\Phi:\prod_{i=1}^{2}\mathcal{P}(\Gamma^{i}_{\R}) \to 2^{\prod_{i=1}^{2}\mathcal{P}(\Gamma^{i}_{\R})}$ as $\Phi=\Psi \circ \Theta$ with $\Theta(P_{\pi,\R}^{1}, P_{\pi,\R}^{2}):=(\Lambda^{1}, \Lambda^{2})$ satisfying the consistency condition \eqref{eq:consistency}, and $\Psi(\Lambda^{1}, \Lambda^{2})=\prod_{i=1}^{2}\text{BR}({P}_{\pi,\R}^{-i}, \Lambda^{1}, \Lambda^{2})$. We have $\mathcal{P}(\Gamma^{i}_{\R})$ is non-empty, compact and convex. Next, we show that the graph $$G=\left\{\left((P_{\pi,\R}^{1}, P_{\pi,\R}^{2}),  \Phi(P_{\pi,\R}^{1}, P_{\pi,\R}^{2})\right)\bigg| (P_{\pi,\R}^{1}, P_{\pi,\R}^{2}) \in \prod_{i=1}^{2}\mathcal{P}(\Gamma^{i}_{\R})\right\}$$ is closed. 

Suppose that sequences $\{P_{\pi,\R}^{1, n}\}_{n}$ and $\{P_{\pi,\R}^{2,n}\}_{n}$ converge to $P_{\pi,\R}^{1, \infty}$ and $P_{\pi,\R}^{2,\infty}$, respectively. Let $P_{\pi,\R}^{i, n}\in \text{BR}({P}_{\pi,\R}^{-i,n}, \Lambda^{1}_{n}, \Lambda^{2}_{n})$ for $i=1,2$ and $\Theta(P_{\pi,\R}^{1,n}, P_{\pi,\R}^{2,n}):=(\Lambda^{1}_{n}, \Lambda^{2}_{n})$. Since the distribution of $u^{i,n}_{\R}$, induced by the policy $P_{\pi,\R}^{i, n}$ converges to the distribution of $u^{i,\infty}_{\R}$, induced by the policy $P_{\pi,\R}^{i, \infty}$, using the fact that $\Theta(P_{\pi,\R}^{1}, P_{\pi,\R}^{2})$ is singleton, we get that $\Theta$ is continuous. Next, we show that $\Psi$ is upper-hemicontinous. We have for $i=1,2$,
\begin{align}
&\inf_{\hat{P}_{\pi,\R}^{i}\in \mathcal{P}(\Gamma^{i}_{\R})} \int c^{i}\left(\omega_{0},u^{i}_{\R},\Xi^{i}\left(\Lambda^{i}_{\infty}(du^{i}_{\R})\right), \Xi^{-i}\left(\Lambda^{-i}_{\infty}(du^{-i}_{\R})\right)\right)\Lambda^{i}_{\infty}(du^{i}_{\R}) \nonumber\\
&\qquad \times \hat{P}_{\pi, \R}^{i}(d{\gamma}^{i}_{\R})P_{\pi,\R}^{-i,\infty}(d{\gamma}^{-i}_{\R})\mathbb{P}^{0}(d\omega_{0})\prod_{s=1}^{2}\hat{f}^{s}\left({y}^{s}_{\R}, \omega_{0}\right)Q^{s}(dy^{s}_{\R})\gamma_{\R}^{s}(du^s_{\R}|y^s_{\R})\nonumber\\
&=\inf_{\hat{P}_{\pi,\R}^{i}\in \mathcal{P}(\Gamma^{i}_{\R})} \lim_{n\to \infty}\int c^{i}\left(\omega_{0},u^{i}_{\R},\Xi^{i}\left(\Lambda^{i}_{n}(du^{i}_{\R})\right), \Xi^{-i}\left(\Lambda^{-i}_{n}(du^{-i}_{\R})\right)\right) \Lambda^{i}_{n}(du^{i}_{\R})\label{eq:pf-the2-36}\\
&\qquad \times \hat{P}_{\pi, \R}^{i}(d{\gamma}^{i}_{\R})P_{\pi,\R}^{-i, n}(d{\gamma}^{-i}_{\R})\mathbb{P}^{0}(d\omega_{0})\prod_{s=1}^{2}\hat{f}^{s}\left({y}^{s}_{\R}, \omega_{0}\right)Q^{s}(dy^{s}_{\R})\gamma_{\R}^{s}(du^s_{\R}|y^s_{\R})\nonumber\\
&\geq \int c^{i}\left(\omega_{0},u^{i}_{\R},\Xi^{i}\left(\Lambda^{i}_{\infty}(du^{i}_{\R})\right), \Xi^{-i}\left(\Lambda^{-i}_{\infty}(du^{-i}_{\R})\right)\right)\Lambda^{i}_{\infty}(du^{i}_{\R}) \label{eq:pf-the2-38}\\
& \qquad \times {P}_{\pi, \R}^{i, \infty}(d{\gamma}^{i}_{\R})P_{\pi,\R}^{-i, \infty}(d{\gamma}^{-i}_{\R})\mathbb{P}^{0}(d\omega_{0})\prod_{s=1}^{2}\hat{f}^{s}\left({y}^{s}_{\R}, \omega_{0}\right)Q^{s}(dy^{s}_{\R})\gamma_{\R}^{s}(du^s_{\R}|y^s_{\R})\nonumber,
\end{align}
where \eqref{eq:pf-the2-36} follows from the generalized dominated convergence theorem for varying measures since the bounded costs $c^{i}$ converge continuously and the actions spaces are compact. Inequality \eqref{eq:pf-the2-38} follows from exchanging limit and infimum and since $P_{\pi,\R}^{i, n}\in \text{BR}({P}_{\pi,\R}^{-i,n}, \Lambda^{1}_{n}, \Lambda^{2}_{n})$, and along the same line as \eqref{eq:pf-the2-36}. Hence, the above chain of inequalities become the chain of equalities, and this implies that $G$ is closed. Moreover, since $\mathcal{P}(\Gamma^{i}_{\R})$ is compact and $\bar{J}^{i}_{\pi,\infty}$ is upper-hemicontinuous, we have that $\text{BR}({P}_{\pi,\R}^{-i}, \Lambda^{1}, \Lambda^{2})$ is non-empty. The Kakutani-Fan-Glicksberg fixed point theorem completes the proof of this step.

\item [{\it Step 2.}] 
It suffices to show that for fixed $P^{-i*}_{\pi,\infty}$ and $\Lambda^{-i}_{\infty}$, the following equality holds:
    \begin{flalign}
    \inf_{P_{\pi}^{i} \in \LPRS^{i}}  J_{\pi, \infty}^{i} (P_{\pi}^{i},{P}_{\pi, \infty}^{-i*})= \inf_{P_{\pi}^{i}\in L^{i}}  J_{\pi, \infty}^{i} (P_{\pi}^{i},{P}_{\pi, \infty}^{-i*}).\label{eq:step2}
    \end{flalign}
    By fixing $P^{-i*}_{\pi,\infty}$, PL$^{i}$ faces a mean-field team problem, and hence,  \eqref{eq:step2} essentially follows from an argument used in \cite[Theorem 2]{SSYdefinetti2020} using \cite[Lemma 2 and Theorem 1]{SSYdefinetti2020}. For completeness, we included some details of the proof.
    
    Let $\LEX^i|_{N}$ be the set of  randomized policies in $\LEX^i$ restricted to their first $N_i$ components. By Lemma \ref{the:0}, we have $\LEX^i|_{N}=\LCOSN^i$, but $\LEX^i|_{N}\not=\LEXN^i$. We have,
    \begin{flalign}
    &\inf_{\hat{P}_{\pi}^{i}\in L^{i}}\limsup_{N_{1},N_{2}\to \infty}\int \frac{1}{N_{i}}\sum_{k=1}^{N_{i}}c^{i}\left(\omega_{0},u^{i}_{k}, \Xi^{1}\left(\frac{1}{N_{1}}\sum_{k=1}^{N_{1}}u^{1}_{k}\right), \Xi^{2}\left( \frac{1}{N_{2}}\sum_{k=1}^{N_{2}}u^{2}_{k}\right)\right) \nonumber\\
    &\qquad \times \hat{P}_{\pi}^{i}(d\underline{\gamma}^{i})P_{\pi,\R}^{-i*,\infty}(d\underline{\gamma}^{-i})\mathbb{P}^{0}(d\omega_{0})\prod_{i=1}^{2}\prod_{k=1}^{N_{i}}\hat{f}^{i}(d{y}^{i}_{k}, \omega_{0})Q^{i}(d{y}^{i}_{k})\gamma_{k}^{i}(du^i_{k}|y^i_{k})\nonumber\\
    &\geq\limsup_{N_{1}, N_{2}\to \infty}\inf_{\hat{P}_{\pi}^{i}\in \LEXN^{i}}\int \frac{1}{N_{i}}\sum_{k=1}^{N_{i}}c^{i}\left(\omega_{0},u^{i}_{k}, \Xi^{1}\left(\frac{1}{N_{1}}\sum_{k=1}^{N_{1}}u^{1}_{k}\right), \Xi^{2}\left( \frac{1}{N_{2}}\sum_{k=1}^{N_{2}}u^{2}_{k}\right)\right) \label{eq:st2-the2-41}\\
    &\qquad \times \hat{P}_{\pi}^{i}(d\underline{\gamma}^{i})P_{\pi,\R}^{-i*,\infty}(d\underline{\gamma}^{-i})\mathbb{P}^{0}(d\omega_{0})\prod_{i=1}^{2}\prod_{k=1}^{N_{i}}\hat{f}^{i}(d{y}^{i}_{k}, \omega_{0})Q^{i}(d{y}^{i}_{k})\gamma_{k}^{i}(du^i_{k}|y^i_{k})\nonumber\\
    &=\limsup_{N_{1},N_{2}\to \infty}\inf_{\hat{P}_{\pi}^{i}\in \LEX^{i}|_{N_i}}\int \frac{1}{N_{i}}\sum_{k=1}^{N_{i}}c^{i}\left(\omega_{0},u^{i}_{k}, \Xi^{1}\left(\frac{1}{N_{1}}\sum_{k=1}^{N_{1}}u^{1}_{k}\right), \Xi^{2}\left( \frac{1}{N_{2}}\sum_{k=1}^{N_{2}}u^{2}_{k}\right)\right) \label{eq:st2-the2-42}\\
    &\qquad \times \hat{P}_{\pi}^{i}(d\underline{\gamma}^{i})P_{\pi,\R}^{-i*,\infty}(d\underline{\gamma}^{-i})\mathbb{P}^{0}(d\omega_{0})\prod_{i=1}^{2}\prod_{k=1}^{N_{i}}\hat{f}^{i}(d{y}^{i}_{k}, \omega_{0})Q^{i}(d{y}^{i}_{k})\gamma_{k}^{i}(du^i_{k}|y^i_{k})\nonumber\\
    &=\limsup_{N_{1}, N_{2}\to \infty}\inf_{\hat{P}_{\pi}^{i}\in \LPRS^{i}|_{N_{i}}}\int \frac{1}{N_{i}}\sum_{k=1}^{N_{i}}c^{i}\left(\omega_{0},u^{i}_{k}, \Xi^{1}\left(\frac{1}{N_{1}}\sum_{k=1}^{N_{1}}u^{1}_{k}\right), \Xi^{2}\left( \frac{1}{N_{2}}\sum_{k=1}^{N_{2}}u^{2}_{k}\right)\right) \label{eq:st2-the2-43}\\
    &\qquad \times \hat{P}_{\pi}^{i}(d\underline{\gamma}^{i})P_{\pi,\R}^{-i*,\infty}(d\underline{\gamma}^{-i})\mathbb{P}^{0}(d\omega_{0})\prod_{i=1}^{2}\prod_{k=1}^{N_{i}}\hat{f}^{i}(d{y}^{i}_{k}, \omega_{0})Q^{i}(d{y}^{i}_{k})\gamma_{k}^{i}(du^i_{k}|y^i_{k})\nonumber\\
    &\geq\inf_{P_{\pi}^{i} \in \LPRS^{i}}  J_{\pi, \infty}^{i} (P_{\pi}^{i},{P}_{\pi, \infty}^{-i*})\label{eq:st2-the2-44}
    \end{flalign}
    where \eqref{eq:st2-the2-41} follows from exchanging limsup and infimum, and step 2 of the proof of Theorem \ref{the:1}. This is because, by fixing the policy of PL$^{-i}$, PL$^i$ faces a team problem with an exchangeable cost function and IS. Hence, without loss of generality, we can restrict our search space to exchangeable policies. Equality \eqref{eq:st2-the2-42} follows from  \cite[Lemma 2]{SSYdefinetti2020} for teams, showing that in the limit, without loss of optimality, the search for (globally) optimal policies of team problems can be restricted to the set of infinitely exchangeable randomized policies. In other words, although $\LEX^i|_{N_i}\subset \LEXN^i$, in the limit, the optimal costs become identical. Equality \eqref{eq:st2-the2-43} follows from the de Finetti representation theorem for infinitely exchangeable policies \cite[Theorem 1]{SSYdefinetti2020}, and linearity of the expected cost in randomized policies $\hat{P}_{\pi}^{i}$ since $\LPRS|_{N_{i}}$ is the set of extreme points of $\LCOS|_{N_{i}}$. Next, we show that \eqref{eq:st2-the2-44} holds.

    This is because, by fixing the policy of PL$^{-i}$, PL$^i$ faces a team problem, and hence, since the cost function is continuous and actions space is compact, by \cite[Theorem 5.2]{yuksel2018general}, a minimizer of the expected cost in \eqref{eq:st2-the2-43} for every finite $N_{1}$ and $N_{2}$ exists. Denote such an optimal policy by $\{\tilde{P}_{\pi, N_{1,2}}^{i}\}_{N_{1,2}}\subseteq \LPRS^{i}|_{N_{i}}$. Since  $\LPRS^{i}|_{N_{i}}$ is compact, there exists a subsequence $\{\tilde{P}_{\pi, n}^{i}\}_{n}\subseteq \LPRS^{i}|_{N_{i}}$ with index $n\in \mathbb{I}\subseteq \N_{i}$ that converges to $\tilde{P}_{\pi, \infty}^{i} \in \LPRS^{i}$. Since $N_{1}$ and $N_{2}$ converge at the same rate, as $n$ drives to infinity, $N_{1}$ and $N_{2}$  drive also to infinity. Define the empirical measures on actions as 
     \begin{flalign}
     \tilde{\Lambda}^{i}_{n}(\cdot):=\frac{1}{n}\sum_{k=1}^{n}\delta_{\{\tilde{u}_{k, n}^{i}\}}(\cdot), \qquad \tilde\Lambda^{i}_{n, \infty}(\cdot):=\frac{1}{n}\sum_{k=1}^{n}\delta_{\{\tilde{u}_{k, \infty}^{i}\}}(\cdot)\label{eq:emp1},
     \end{flalign}
     where actions $\tilde{u}_{k, n}^{i}$ and$\tilde{u}_{k, \infty}^{i}$ are induced by randomized policies $\tilde{P}_{\pi, n}^{i}$ and $\tilde{P}_{\pi, \infty}^{i}$, respectively.
     Similarly, we define $\Lambda^{-i}_{n, \infty}$ where the actions ${u}_{k, \infty}^{-i}$ are induced by randomized policies ${P}_{\pi, \R}^{-i*, \infty}$. By Assumption \ref{assump:oc} and under the change of measure in \eqref{eq:abscon} (see Assumption \ref{assump:ind}), the observations $(y^{i}_{k})_{k\in \N_{i}}$ are i.i.d. Since randomized policies $\tilde{P}_{\pi, n}^{i}$ and $\tilde{P}_{\pi, \infty}^{i}$ belong to $\LPRS^{i}|_{n}$ and $\LPRS^{i}$, respectively, we can conclude that $\tilde{u}_{k, \infty}^{i}$, $\tilde{u}_{k, \infty}^{i}$, and ${u}_{k, \infty}^{-i}$ are i.i.d. Hence, by the strong law of large numbers, along the same lines as in the proof of \cite[Theorem 3.2]{sanjari2019optimal} (see \cite[Eqs (3.6)--(3.10)]{sanjari2019optimal}),  sequences $\{\tilde\Lambda^{i}_{n}\}_{n}$, and $\{\tilde{\Lambda}^{i}_{n,\infty}\}_{n}$ converge weakly to the same limit $\tilde\Lambda^{i}_{\infty}:=\mathcal{L}(u^{i*}_{\R, \infty})$.   Also, $\Lambda^{-i}_{n, \infty}$ converges weakly to $\Lambda^{-i}_{\infty}:=\mathcal{L}(u^{-i*}_{\R, \infty})$.
    Using Assumption \ref{assump:cc}, and the generalized dominated convergence theorem, we get \eqref{eq:st2-the2-44}, and complete the proof.

\end{itemize}

\subsection{Proof of Theorem \ref{the:5}.}\label{app:E}
By Theorem \ref{the:2}, there exists a symmetric NE $(P_{\pi}^{1*},P_{\pi}^{2*})\in \LPRS^{1}\times \LPRS^{2}$ for \PIN. Proceeding along the same lines as \eqref{eq:st2-the2-43}--\eqref{eq:st2-the2-44}, we get
\begin{align}
    &J_{\pi, \infty}^{i}(P_{\pi}^{1*},P_{\pi}^{2*})=\limsup_{N_{i}\to \infty} \inf_{P_{\pi, N}^{i}\in \LPRS^{i}|_{N_{i}}} J_{\pi, N}^{i} (P_{\pi}^{i},{P}_{\pi}^{-i*}),\quad \forall i=1,2.\label{eq:56-the3a}
\end{align}
Following \cite[Theorem 5.2]{yuksel2018general}, there exists a sequence of optimal policies as a best response to ${P}_{\pi}^{-i*}$, denoted by $\{\tilde{P}_{\pi, N_i}^{i}\}_{N_{i}}$. Again proceeding along the same lines as \eqref{eq:st2-the2-43}--\eqref{eq:st2-the2-44} with replacing the limsup with liminf, we get for all $i=1,2$,
\begin{align}
    \liminf_{N_{i}\to \infty} \inf_{P_{\pi, N}^{i}\in \LPRS^{i}|_{N_{i}}} J_{\pi, N}^{i} (P_{\pi}^{i},{P}_{\pi}^{-i*})
    &\geq\inf_{P_{\pi}^{i} \in \LPRS^{i}}  J_{\pi, \infty}^{i} (P_{\pi}^{i},{P}_{\pi}^{-i*})\label{eq:a56-the3}\\
    &=J_{\pi, \infty}^{i}(P_{\pi}^{1*},P_{\pi}^{2*})\nonumber,
\end{align}
where \eqref{eq:a56-the3} follows from an argument similar to that used in \eqref{eq:st2-the2-43}--\eqref{eq:st2-the2-44}. This together with \eqref{eq:56-the3a}, implies that
\begin{align}
    &J_{\pi, \infty}^{i}(P_{\pi}^{1*},P_{\pi}^{2*})=\lim_{N_{i}\to \infty} \inf_{P_{\pi, N}^{i}\in \LPRS^{i}|_{N_{i}}} J_{\pi, N}^{i} (P_{\pi}^{i},{P}_{\pi}^{-i*}),\quad \forall i=1,2.\label{eq:56-the3}
\end{align}
By the strong law of large numbers since $(P_{\pi}^{1*},P_{\pi}^{2*})\in \LPRS^{1}\times \LPRS^{2}$, we get
\begin{align}
    &J_{\pi, \infty}^{i}(P_{\pi}^{1*},P_{\pi}^{2*})=\lim_{N_{i}\to \infty} J_{\pi, N}^{i}(P_{\pi}^{1*},P_{\pi}^{2*}),\quad \forall i=1,2.\label{eq:57-the3}
\end{align}
Equality \eqref{eq:56-the3} together with \eqref{eq:57-the3} implies that there exists $\hat\epsilon_{N_{1:2}}>0$, converging to zero as  $N_{1},N_{2}\to \infty$ such that
\begin{align}
    &\left|J_{\pi, N}^{i}(P_{\pi}^{i*},P_{\pi}^{-i*})- \inf_{P_{\pi, N}^{i}\in \LPRS^{i}|_{N_{i}}} J_{\pi, N}^{i} (P_{\pi}^{i},{P}_{\pi}^{-i*})\right|\leq \hat\epsilon_{N_{1:2}},\quad \forall i=1,2.\label{eq:the3-54}
\end{align}
This implies that $(P_{\pi}^{1*},P_{\pi}^{2*})$ constitutes an $\hat\epsilon_{N_{1:2}}$-NE for \PN\ among symmetric policies $\LPRS^{1}|_{N_{1}}\times \LPRS^{2}|_{N_{2}}$. Proceeding along the same lines as \eqref{eq:st2-the2-41}--\eqref{eq:st2-the2-44}, we get 
\begin{flalign*}
\limsup_{N_{i}\to \infty} \inf_{P_{\pi, N}^{i}\in L^{i}_{N}}  J_{\pi, N}^{i} (P_{\pi}^{i},{P}_{\pi}^{-i*})=\limsup_{N_{i}\to \infty} \inf_{P_{\pi, N}^{i}\in \LPRS^{i}|_{N_{i}}}  J_{\pi, N}^{i} (P_{\pi}^{i},{P}_{\pi}^{-i*})\quad \forall i=1,2.
\end{flalign*}
Hence, using \eqref{eq:the3-54}, there exists $\epsilon_{N_{1:2}}>0$, converging to zero as  $N_{1},N_{2}\to \infty$ such that
\begin{align}
    &\left|J_{\pi, N}^{i}(P_{\pi}^{i*},P_{\pi}^{-i*})- \inf_{P_{\pi, N}^{i}\in L^{i}_{N}} J_{\pi, N}^{i} (P_{\pi}^{i},{P}_{\pi}^{-i*})\right|\leq \epsilon_{N_{1:2}},\quad \forall i=1,2.\label{eq:the3-56}
\end{align} 
This implies that $(P_{\pi}^{1*},P_{\pi}^{2*})$ constitutes an $\epsilon_{N_{1:2}}$-NE for \PN\ among all policies in $L^{1}_{N}\times L^{2}_{N}$, and this completes the proof.

\subsection{Proof of Theorem \ref{the:3}.}\label{app:C}

The proof is similar to  that of Theorem \ref{the:1}. We first use the Kakutani-Fan-Glicksberg fixed point theorem to show that there exists a NE $(P_{\pi}^{1*}, P_{\pi}^{2*})$ among all policies in $\LEXN^{1} \times \LEXN^{2}$. This follows from an identical argument as that in the step 1 of the proof of Theorem \ref{the:1} using  continuity of the costs and observations in actions in Assumptions \ref{assump:sxcdy} (ii) and (iii).  Then, we can show that an equality \eqref{eq:e-cost} adapted for \PNT\ holds using Assumptions \ref{assump:sxcdy} (iv) and (v), and an argument similar to \cite[Lemma 1]{SSYdefinetti2020} (see \cite[Lemma 3]{SSYdefinetti2020} for details of this step). For brevity, the details are not included.

\subsection{Proof of Theorem \ref{the:4}.}\label{app:D}
The proof proceeds along the same lines as those in the proof of Theorem \ref{the:2}. In the following, we only present some additional technical arguments required for the dynamic setting.  We first restrict the search in finding a NE for \PIT\ to randomized policies belonging to $\LPRS^{1} \times \LPRS^{2}$. By the strong law of large numbers and under Assumptions \ref{assump:c}, \ref{assump:2}, \ref{assump:ind1} (i) and \ref{assump:ind2} hold, we can rewrite the expected costs $J_{\pi, \infty}^{i}(P_{\pi}^{1},P_{\pi}^{2})$ in \eqref{eq:infiniterandom} for any $(P_{\pi}^{1},P_{\pi}^{2})$ in $\LPRS^{1} \times \LPRS^{2}$ as $\bar{J}_{\pi, \infty}^{i,T}(P_{\pi, \R}^{i},P_{\pi, \R}^{-i}, {\Lambda}^{i}, {\Lambda}^{-i})$ with
    $P_{\pi}^{i}(\pmb{\underline\gamma^{i}} \in \cdot)=\prod_{k\in \mathbb{N}}{P}_{\pi, \R}^{i}(\pmb{\gamma^{i}_{k}} \in \cdot)$,
for some ${P}_{\pi, \R}^{i}(\pmb{\gamma^{i}_{k} }\in \cdot)$ belonging to the set $\mathcal{P}(\prod_{t=0}^{T-1}\Gamma^{i}_{\R,t})$, and 
\begin{align}
\Lambda^{i}(\cdot \times \prod_{t=0}^{T-1} \mathbb{X}^{i} \times \mathbb{Y}^{i} \times \mathbb{Z}^{i})=\mathcal{L}(\pmb{u^{i}_{\R}}), \qquad \Lambda^{i}(\cdot \times \prod_{t=0}^{T-1} \mathbb{U}^{i} \times \mathbb{Y}^{i} \times \mathbb{Z}^{i})=\mathcal{L}(\pmb{x^{i}_{\R}}).\label{eq:consistency-dyn}
\end{align}
Similarly, we can rewrite the expected cost under Assumption \ref{assump:ind1} (ii). Similar to the proof of Theorem \ref{the:2}, we define $\Phi:\prod_{i=1}^{2}\mathcal{P}(\prod_{t=0}^{T-1}\Gamma^{i}_{\R,t}) \to 2^{\prod_{i=1}^{2}\mathcal{P}(\prod_{t=0}^{T-1}\Gamma^{i}_{\R,t})}$ as the best response correspondence $\Phi=\Psi \circ \Theta$ with $\Theta(P_{\pi,\R}^{1}, P_{\pi,\R}^{2}):=(\Lambda^{1}, \Lambda^{2})$ satisfying the consistency condition \eqref{eq:consistency-dyn}, and $\Psi(\Lambda^{1}, \Lambda^{2})=\prod_{i=1}^{2}\text{BR}({P}_{\pi,\R}^{-i}, \Lambda^{1}, \Lambda^{2})$. 
The set $\mathcal{P}(\prod_{t=0}^{T-1}\Gamma^{i}_{\R,t})$ is non-empty, compact, and convex. Next, we show that the graph $G$ is closed.

Suppose that sequences $\{P_{\pi,\R}^{1, n}\}_{n}$ and $\{P_{\pi,\R}^{2,n}\}_{n}$ converge to $P_{\pi,\R}^{1, \infty}$ and $P_{\pi,\R}^{2,\infty}$, respectively. Let $P_{\pi,\R}^{i, n}\in \text{BR}({P}_{\pi,\R}^{-i,n}, \Lambda^{1}_{n}, \Lambda^{2}_{n})$ for $i=1,2$ and $\Theta(P_{\pi,\R}^{1,n}, P_{\pi,\R}^{2,n}):=(\Lambda^{1}_{n}, \Lambda^{2}_{n})$. Since the distribution of $u^{i,n}_{\R,0}$, induced by the policy $P_{\pi,\R}^{i, n}$ converges to the distribution of $u^{i,\infty}_{\R.0}$, induced by the policy $P_{\pi,\R}^{i, \infty}$, the marginal of $\Lambda^{i}_{n}$ on $u_{\R,0}^{i,n}$ converges to the marginal of $\Lambda^{i}_{\infty}$ on $u_{\R,0}^{i,\infty}$. By continuity of $f_{t}^{i}$, $\Xi^{i}_{x}$, and $\Xi^{i}_{u}$, we get that $x^{i,n}_{\R,1}$ converges weakly to  $x^{i,\infty}_{\R.1}$. Hence, marginals of $\Lambda^{i}_{n}$ on $x_{\R,1}^{i,n}$ converges to marginals of $\Lambda^{i}_{\infty}$ on $x_{\R,1}^{i,\infty}$. Hence, we can recursively show that $\Lambda_{n}^{i}$ converges weakly to $\Lambda_{\infty}^{i}$ for $i=1,2$ (and $\mathbb{P}^{0}$-a.s.). This leads to the continuity of $\Theta$. Along the same lines as \eqref{eq:pf-the2-36}--\eqref{eq:pf-the2-38}, using the generalized dominated convergence theorem, it can be shown that $\Psi$ is upper-hemicontinous. Hence, by the Kakutani-Fan-Glicksberg fixed point theorem, there exists a policy profile $({P}_{\pi, \infty}^{1*}, {P}_{\pi, \infty}^{2*})$ that constitutes a NE for \PIT, among all randomized policies in $\LPRS^{1} \times \LPRS^{2}$. By fixing $P^{-i*}_{\pi,\infty}$ and $\Lambda^{-i}_{\infty}$, PL$^{i}$ faces a team problem, and hence, \cite[Theorem 4]{SSYdefinetti2020} implies that \eqref{eq:step2} holds. Hence,  $({P}_{\pi, \infty}^{1*}, {P}_{\pi, \infty}^{2*})$ constitutes a NE for \PIT, among all randomized policies in $L^{1} \times L^{2}$, and completes the proof.

\bibliographystyle{plain}

\end{document}